\documentclass{article}

\usepackage{amsmath,amsthm,amsfonts,amssymb,fullpage}
\usepackage[utf8]{inputenc}

\usepackage{xcolor}
\definecolor{grau}{rgb}{0.1,0.1,0.1}
\usepackage[colorlinks=true, urlbordercolor=grau, linkcolor=grau, citecolor=grau, urlcolor=grau]{hyperref}
\usepackage{accents}
\usepackage{booktabs}

\newtheorem{theorem}{Theorem}[section]

\newtheorem{corollary}[theorem]{Corollary}
\newtheorem{lemma}[theorem]{Lemma}

\newtheorem{proposition}[theorem]{Proposition}

\theoremstyle{definition}

\newtheorem{remark}[theorem]{Remark}

\makeatletter
\renewenvironment{proof}[1][\proofname]{%
   \par\pushQED{\qed}\normalfont%
   \topsep6\p@\@plus6\p@\relax
   \trivlist\item[\hskip\labelsep\bfseries#1\@addpunct{.}]%
   \ignorespaces
}{%
   \popQED\endtrivlist\@endpefalse
}
\makeatother

\newcommand\C{{\mathbb C}}
\newcommand\F{{\mathbb F}}
\newcommand\Q{{\mathbb Q}}
\newcommand\R{{\mathbb R}}
\newcommand\Z{{\mathbb Z}}

\newcommand\GL{{\mathop{\mathrm{GL}}}}
\newcommand\SL{{\mathop{\mathrm{SL}}}}

\newcommand\gal{{\mathop{\mathrm{Gal}}}}
\newcommand\height{{{\mathrm{h}}}}
\newcommand\im{{\mathop{\mathrm{Im}}}}
\newcommand\id{{\mathop{\mathrm{id}}}}
\newcommand\ns{{\mathop{\mathrm{ns}}}}
\newcommand\ord{{\mathop{\mathrm{Ord}}}}
\newcommand\re{{\mathop{\mathrm{Re}}}}
\newcommand\spl{{\mathop{\mathrm{sp}}}}

\newcommand\CC{{\mathcal C}}
\newcommand\HH{{\mathcal H}}

\newcommand\NN{{\mathcal N}}
\newcommand\OO{{\mathcal O}}

\newcommand\bfa{{\mathbf a}}
\newcommand\bfb{{\mathbf b}}

\newcommand\ph\varphi
\newcommand\eps\varepsilon

\newcommand\gerk{{\mathfrak k}}
\newcommand\gerp{{\mathfrak p}}

\newcommand\tila{{\tilde{a}}}
\newcommand\tilb{{\tilde{b}}}
\newcommand\tild{{\tilde{d}}}
\newcommand\tilsigma{{\tilde\sigma}}
\newcommand\tilbfa{{\tilde\bfa}}
\newcommand\tilbfb{{\tilde\bfb}}

\newcommand\hatmho{{\widehat\mho}}

\newcommand{\nmax}{{n_0}}

\newcommand\ocalF{{\accentset{~\circ}{\calF}}}
\newcommand\oOmega{{\accentset{\circ}{\Omega}}}

\newcommand{\SAGE}{\textnormal{SageMath}}
\newcommand{\PARIGP}{\textsc{Pari/GP}}
\newcommand{\RR}{\mathbb{R}} %reals
\newcommand{\QQ}{\mathbb{Q}} %rationals
\newcommand{\ZZ}{\mathbb{Z}} %integers
\DeclareMathOperator*{\wo}{\backslash}
\newcommand{\conv}{\textnormal{conv}} %convex hull
\newcommand{\ceil}[1]{\left\lceil{#1}\right\rceil} %rounding up
\newcommand{\st}{\ |\ } %such that '|'
\newcommand{\wt}[1]{\widetilde{#1}}
\newcommand{\calF}{{\mathcal{F}}}

\usepackage{marginnote}

\newcommand{\addresssize}{\small}

\numberwithin{equation}{section}

\title{Computing integral points on $X_\ns^+(p)$}

\author{Aurélien Bajolet, Yuri Bilu\footnote{Supported by the \textsl{Agence National de la Recherche} project ``Hamot'' (ANR 2010 BLAN-0115-01) and by the ALGANT scholarship program, and by the SPARC Project P445 (India).}, \setcounter{footnote}{6} Benjamin Matschke\footnote{Supported by the Max Planck Institute for Mathematics Bonn, by IdEx project DiGeMANT, by Koç University, and by Simons Foundation grant \#{}550023.}}

1

\makeatletter
\renewcommand*\l@section[2]{%
  \ifnum \c@tocdepth >\z@
    \addpenalty\@secpenalty
    \addvspace{0.3em \@plus\p@}%
    \setlength\@tempdima{1.5em}%
    \begingroup
      \parindent \z@ \rightskip \@pnumwidth
      \parfillskip -\@pnumwidth
      \leavevmode %\bfseries
      \advance\leftskip\@tempdima
      \hskip -\leftskip
      #1\nobreak\hfil \nobreak\hb@xt@\@pnumwidth{\hss #2}\par
    \endgroup
  \fi}
\makeatother

\begin{document}

\hfuzz 4pt

\maketitle

\vspace{0.5em}

\begin{flushright}
\textit{To the memory of Alan Baker}
\end{flushright}

\vspace{0.5em}

\begin{abstract}
We develop a general method for computing integral points on modular curves, based on Baker's inequality. As an illustration, we show that for ${11\le p<101}$, the only integral points  on the curve $X_\ns^+(p)$  are the CM points. 
\end{abstract}

{\footnotesize
\tableofcontents}

\section{Introduction}
In his celebrated work of 1978  Mazur~\cite{Ma78} completely described the possible  rational points on the modular curves $X_0(p)$, where~$p$ is a prime number. In particular, he showed that the set $X_0(p)(\Q)$ consists only of the cusps if ${p>163}$, and of the cusps and the CM-points if ${37<p\le 163}$. 

The curve $X_0(p)$ is associated to the Borel subgroup of $\GL_2(\F_p)$.
It is natural to ask the same question on the modular curves associated to two other important maximal subgroups of $\GL_2(\F_p)$, the normalizers of a split Cartan or a non-split Cartan subgroup. See \cite[Appendix~A.5]{Se97} or \cite[Section~2]{Ba10}, where all the necessary definitions are given.
We shall denote these curves $X_\spl^+(p)$ and $X_\ns^+(p)$, respectively\footnote{In~\cite{Se97} the curves $X_\spl^+(N)$ and  $X_\ns^+(N)$ are defined  but for arbitrary levels~$N$, but in this article we restrict to prime levels only.}.
This problem is not only interesting by itself, but is also motivated by applications; for instance, Serre's uniformity problem about Galois representations~\cite{BP11} would be solved if one could show that for large~$p$ the sets $X_\spl^+(p)(\Q)$ and $X_\ns^+(p)(\Q)$ consist only of the cusps and the CM-points (points corresponding to elliptic curves with complex multiplication).
For the convenience of the reader, we reproduce the full list of the 13  rational CM $j$-invariants in Table~\ref{tacm}. 

\begin{table}[t]
{\scriptsize
$$
\begin{gathered}
\begin{array}{c|ccccccccccccc}
df^2&
-3& 
-3\cdot 2^2 &
-3\cdot3^2&
-4&
-4\cdot2^2&
-7&
-7\cdot2^2&
-8\\
j&
0 &
2^43^35^3& 
-2^{15}3^15^3& 
2^63^3& 
2^33^311^3& 
-3^35^3& 
3^35^317^3&
2^65^3
\end{array}\\
\begin{array}{c|ccccccccccccc}
df^2 &
-11 &
-19 &
-43 &
-67 &
-163\\
j& 
-2^{15}& 
-2^{15}3^3 &
-2^{18}3^35^3 &
-2^{15}3^35^311^3 &
-2^{18}3^35^323^329^3 
\end{array}
\end{gathered}
$$
\caption{Rational CM $j$-invariants together with the discriminant of the CM order.}
\label{tacm} 
}
\end{table}

Rational points on the curves $X_\spl^+(p)$ were determined in \cite{BP11a,BPR13} for all ${p\ne 13}$; in particular, it is shown in~\cite{BPR13} that for ${p\ge 17}$ the set $X_\spl^+(p)(\Q)$  consists only of the cusps and the CM-points. The case ${p=13}$ was resolved  in a recent breakthrough  by Balakrishnan et al~\cite{BDMTV19}, who computed all rational points on  $X_\spl^+(13)$, using Kim's ``quadratic Chabauty'' method.

Unfortunately, the methods of \cite{BP11a,BPR13} completely fail for the curve $X_\ns^+(p)$. To the best of our knowledge,  the set $X_\ns^+(p)(\Q)$ is not known for any prime ${p\ge 17}$.

More is known about \textsl{integral} points on the curves $X_\ns^+(p)$, that is, points ${P\in X_\ns^+(p)(\Q)}$ such that ${j(P)\in \Z}$, where~$j$ is the modular invariant. It is easy to see that for ${p \le 5}$ this set is infinite. Kenku~\cite{Ke85} determined the integral points\footnote{Kenku's list has two typos: instead of ${j = -2^{15}3^35^311^3}$ he writes its negative; and instead of ${j =  2^{15}  7^5}$ he writes ${j=7^5 2^5}$.} on the curve $X_\ns^+(7)$; in fact, he found the $7$-integral points, that is, the points $P\in X_\ns^+(p)(\Q)$ such that the denominator of $j(P)$ is a power of~$7$. He used in an essential way the fact that the curve is of genus~$0$.

More recently, Schoof and Tzanakis~\cite{ST12} determined the integral points  on $X_\ns^+(11)$, using the fact that this curve is of genus~$1$. They showed that the only integral  points on this curve are the CM-points.  See also~\cite{CC04}.

The methods of~\cite{Ke85,ST12} are quite ad hoc and do not extend to other levels. 

Since the curves  $X_\spl^+(13)$ and  $X_\ns^+(13)$ are known to be isomorphic over~$\Q$ (see Baran~\cite{Ba14}), the already mentioned result of 
Balakrishnan et al~\cite{BDMTV19} computes \textsl{rational}  points on $X_\ns^+(13)$ as well. The approach of this spectacular work has some potential of extending to higher levels, but this would require substantial new ideas.

We may also mention that integral points on the curve $X_\ns^+(N)$ of certain composite levels~$N$ were determined much earlier by Heegner and Siegel~\cite{He52,Si68} in the context of the Class Number~$1$ problem;  see \cite[Appendix~A.5]{Se97} for more details. More recently, composite levels were examined by Baran~\cite{Ba09,Ba10}. None of these methods seems to extend to higher prime levels either.

%One can also mention the work of Baran~\cite{Ba10}, who considered some  composite levels as well. 

%It is shown in \cite{Bi95,Bi02,BI11} that, using Baker's method, one can, in principle, determine integral points on ``many'' modular curves; for instance, on those having at least~$3$ cusps  \cite{Bi95,Bi02} or those of prime power level (with a few exceptions)~\cite{BI11}. Moreover, 

It was observed in \cite{Bi95,Bi02} that heights of integral points on a modular curve with~$3$ or more cusps can  be effectively bounded using Baker's method; moreover, this is also true for $S$-integral points defined over an arbitrary number field. 

The modular curve $X_\ns^+(p)$ has $\ceil{(p-1)/2}$ cusps, so the observation above applies as soon as ${p\ge 7}$. Indeed, Bajolet and Sha~\cite{BS13} obtained a fully explicit upper bound for the size of an integral point~$P$ on $X_\ns^+(p)$ for an arbitrary prime ${p\ge 7}$. They showed that in general,
\begin{equation}
\label{ebsh}
\log|j(P)| <41993\cdot 13^{p} \cdot p^{2p+7.5}(\log{p})^{2},
\end{equation}
and this bound can be substantially refined if ${p-1}$ is divisible by a small odd prime or by~$8$, see Theorem~\ref{thbasha} below.  
Sha~\cite{Sh13,Sh14}  extended this result of~\cite{BS13} to $S$-integral points on rather general modular curves over arbitrary number fields, giving an explicit version of the ``effective Siegel theorem for modular curves''~\cite{Bi02,BI11}. See also the recent work of Cai~\cite{Ca19}.

Using numerical bound~\eqref{ebsh}, one can in principle enumerate all integral points on $X_\ns^+(p)$. However, this bound is too huge to perform this enumeration in reasonable time. 

It turns out that the huge bound can be reduced using the numerical Diophantine approximation techniques, which go back to the work of Baker and Davenport~\cite{BD69}. The idea of Baker and Davenport was elaborated in~\cite{BH96,BH98,BH99,BHV01,Ha00,PS87,TW89} in the context of the Diophantine equations of Thue and of related types, providing practical methods for solving these equations.

In the present article we adapt these techniques to modular curves and  develop an algorithm for finding integral points on  the modular curve $X_\ns^+(p)$, where ${p\ge 7}$ is an arbitrary prime number. Having implemented our algorithm, we prove the following.

\begin{theorem}
\label{thmMain}
Let~$p$ be a prime number, ${11\le p \le 97}$, and let ${P\in X_\ns^+(p)(\Q)}$ be such that ${j(P)\in \Z}$. Then~$P$ is a CM point, that is, $j(P)$ is one of the 13 numbers displayed in the second line of Table~\ref{tacm}. 
\end{theorem}

One may conjecture that for any prime ${p\ge 11}$ the only integral points on $X_\ns^+(p)$ are the CM-points. 

It might be useful to recall the description of integral CM-points on $X_\ns^+(p)$. Let~$D$ be a negative quadratic discriminant with ${h(D)=1}$ (that is, one of the 13 numbers in the upper row of Table~\ref{tacm}) and let~$j_D$ be the corresponding $j$-invariant. Then, for ${p\ge3}$, there is a point ${P\in X_\ns^+(p)(\Z)}$ with ${j(P)=j_D}$ if and only if  ${(D/p)=-1}$.

As compared to the previous work, our article has two new features. The first one is using ``economical'' modular units made of Siegel functions. Applying explicit modular units constructed as multiplicative combinations of Siegel functions to Diophantine problems involving modular curves goes back to the work of Kubert and Lang in 70s \cite{KL75}, \cite[Chapter~8]{KL81}. More recently, this was succefully implemented in \cite{BS13,BP11,BP11a,Sh14a,Sh14}. However, the units used therein, while perfect for theoretical bounds, are no longer suitable for numerical purposes. We replace them with ``economical units'', constructed in Sections~\ref{ssiegel}--\ref{scpu}.

Another new element, as compared to \cite{BH96,BH98,BH99,BHV01,Ha00,PS87,TW89},  is related to the fact that, in our case, the Baker-Davenport method alone is not sufficient to eliminate the false positives. We show that potential integral points correspond to lattice points on a certain analytic curve. To rule them out, we cover this curve with ellipsoids and compute integral points in them using the Fincke-Pohst algorithm~\cite{FP85}. See Sections~\ref{secOutlineOfAlgorithm} and~\ref{secEllipsoidsSieve} for more details.

Note that our method   does not use Jacobian embedding, applying thereby to curves of very high genus: for instance, the genus of   $X_\ns^+(97)$ is $353$. Theoretically, our method works for any prime. We did not go beyond ${p=97}$ because for big~$p$ computations become too long even for modern computers. See more on this in Section~\ref{secRunningTime}. 

One may wonder whether our method can be adapted to more general set up, like computing integral or even $S$-integral points over number fields on arbitrary modular curves with~$3$ or more cusps. Again, theoretically this is possible, but practical implementation may require calculations which are too hard for modern computers. A realistic task is, perhaps, computing rational points on $X_\ns^+(p)$ whose denominator is a power of~$p$, like Kenku~\cite{Ke85} did for ${p=7}$. 

Our source code and data is available at \url{https://github.com/bmatschke/x-nonsplit-plus/}.

\paragraph{Acknowledgments}
Yuri Bilu was supported by the \textsl{Agence National de la Recherche} project ``Hamot'' (ANR 2010 BLAN-0115-01), by the \textsc{Algant} scholarship program and by the SPARC Project P445 (India). 
Benjamin Matschke was supported by the Max Planck Institute for Mathematics Bonn, by the project DiGeMANT funded by the Initiative d'excellence de l'Universit\'e de Bordeaux (IdEx), by Koç University, and by Simons Foundation grant \#{}550023.
We thank Julia Baoulina, Denis Benois, Andreas Enge, Elisa Gorla, Pierre Parent  and Sha Min for useful discussions and suggestions.
We thank the referees for their detailed reviews and valuable comments, that helped us to improve the presentation.

Our algorithms are implemented in the computer algebra system  \SAGE~\cite{sage}, which in turn relies on methods from~\PARIGP~\cite{pari}. 
The final computations were carried out on the PlaFRIM computer cluster, supported by Inria, CNRS (LABRI and IMB), Université de Bordeaux, Bordeaux INP and Conseil Régional d’Aquitaine (see \url{https://www.plafrim.fr/}).

\subsection{Plan of the article}
In Section~\ref{sxg} we recall basic definitions about modular curves. In particular, we review the notions of the \textsl{nearest cusp} and the \textsl{$q$-parameter} at a given cusp, a basic tool in the calculus on modular curves.

In Section~\ref{sbakergen} we give a general informal overview on how Baker's method applies to modular curves, highlighting both theoretical and numerical aspects.  

In Sections~\ref{ssiegel} and~\ref{squad} we revise the theory of  modular units,  an indispensable tool in the Diophantine analysis of modular curves. In Section~\ref{scpu} we apply this general theory in the special case of the curve $X_\ns^+(p)$, constructing especially ``economical'' units on this curve. 

In Section~\ref{sprinrel} we evaluate the unit constructed in Section~\ref{scpu} at an integral point~$P$, and express the value as a multiplicative combination of certain algebraic numbers: ${U(P)=\eta_0^{b_0}\eta_1^{b_1}\cdots\eta_r^{b_r}}$. We then express the exponents~$b_k$ in terms of the Galois conjugates of~$U$ and also in terms of the $q$-parameter of~$P$. These expressions, while pretty trivial, will play a fundamental role in the remaining part of the article.

In Section~\ref{secOutlineOfAlgorithm} we outline the algorithm that finds all integral points based on these expressions of $b_k$ using that all $b_k$ are integers. The remaining sections present the different parts of the algorithm in detail.

In Section~\ref{supperboundred} we recall ``Baker's bound'',   a huge explicit upper bound for the $j$-invariants of integral points~$P$ on $X_\ns^+(p)$, obtained in~\cite{BS13} using Baker's method. This implies a very tiny lower bound for the $q$-parameter of~$P$,  and we show how it can be drastically improved in practical situations, using the reduction technique introduced by Baker and Davenport. This way we obtain a more reasonable lower bound for the absolute value of the $q$-parameter, which is still insufficient to list efficiently all integral points just by exhaustive search. 

Therefore in Section~\ref{secEllipsoidsSieve} we present an algorithmic sieve that further reduces this set of possible values for $j(P)$ considerably. It can be seen as a much more detailed elaboration of the previous reduction step, in which some candidates for $j(P)$ may remain, and they may indeed come from integral points.
In Section~\ref{ssmissing} we deal with the possible values of  $j(P)$  left after the sieving. An overview of running times is given in Section~\ref{secRunningTime}.

\subsection{Notation and conventions}
\label{ssconv}

\paragraph{The logarithm.}
Unless the contrary is stated explicitly, for the complex logarithm we choose the branch satisfying 
$$
-\pi <\im\log z \le \pi \qquad (z\in \C^\times). 
$$
Note that, with this definition, we do not always have the equality ${\log(zw)=\log z+\log w}$, but always have the inequality
$$
|\log(zw)|\le |\log z|+|\log w|. 
$$

\paragraph{Modular functions.}
Throughout the article, the letter~$j$ may have four different meanings, sometimes in the same equation, like in~\eqref{ejptau} and~\eqref{ejpq}: the modular invariant $j(\tau)$ on the Poincaré upper halfplane~$\HH$; the modular invariant $j(E)$ of an elliptic curve~$E$; the ``modular invariant'' rational function on a modular curve; the sum of the familiar series ${j(q)=q^{-1}+ 744+ 196884q+\ldots}$. It should be always clear from the context which meaning of~$j$ is used. A similar convention applies to other modular functions as well.

\paragraph{The $O_1(\cdot)$ notation.}
We shall use the notation $O_1(\cdot)$, which is a quantitative  analogue of the familiar $O(\cdot)$. Precisely, ${A=O_1(B)}$ means that ${|A|\le B}$.

\section{Modular curves, nearest  cusps and \texorpdfstring{$q$}q-parameters}
\label{sxg}
Let~$N$ be a positive integer. The modular curve $X(N)$ has a geometrically irreducible model over the cyclotomic field $\Q(\zeta_N)$, and the Galois group 
${\gal\bigl(\Q(\zeta_N)(X(N))\big/\Q(j)\bigr)}$ is canonically isomorphic  to the quotient $\GL_2(\Z/N\Z)/\{\pm1\}$, with $\SL_2(\Z/N\Z)/\{\pm1\}$ being the group ${\gal\bigl(\Q(\zeta_N)(X(N))\big/\Q(\zeta_N,j)\bigr)}$, see \cite[Chapter~6]{La87} or \cite[Sections~7.5 and~7.6]{DS05}. We write the Galois action of $\GL_2(\Z/N\Z)$ on the field $\Q(\zeta_N)\bigl(X(N)\bigr)$ exponentially. In the following proposition we collect the properties of this action.

\begin{proposition}
\label{pgalomod}
\begin{enumerate}
\item
\label{ifsig}
For ${u\in \Q(\zeta_N)\bigl(X(N)\bigr)}$ and ${\sigma \in \SL_2(\Z/N\Z)}$ we have 
\begin{equation*}
u^\sigma= u\circ\tilsigma,
\end{equation*}
where on the right we view~$u$ as a $\Gamma(N)$-automorphic function on the extended Poincar\'e plane $\bar\HH$, and~$\tilsigma$ is a lifting of~$\sigma$ to ${\Gamma(1)=\SL_2(\Z)}$. Clearly, the result is independent of the choice of the lifting. 

\item
\label{idetsigma}
For ${\sigma \in \GL_2(\Z/N\Z)}$ we have 
\begin{equation}
\zeta_N^\sigma=\zeta_N^{\det\sigma}.
\end{equation}

\item
\label{igalqexp}
Recall that ${u\in \Q(\zeta_N)\bigl(X(N)\bigr)}$ has a ``$q$-expansion''
$$
u=\sum_{k=k_0}^\infty a_kq^{k/N}\in\Q(\zeta_N)((q^{1/N})). 
$$
Then for ${\sigma=(\begin{smallmatrix}1&0\\0&d\end{smallmatrix})}$ the $q$-expansion of~$u^\sigma$ is 
\begin{equation*}
u^\sigma=\sum_{k=k_0}^\infty a_k^\sigma q^{k/N}. 
\end{equation*}

\end{enumerate}

\end{proposition}

\begin{proof}
For item~\ref{ifsig} see \cite[bottome of page~280]{DS05}. Item~\ref{idetsigma} is \cite[Lemma~7.6.1]{DS05}; see also \cite[Theorem~3 on page 66]{La87}. 

Item~\ref{igalqexp} is enough to verify for the functions $f_{r,s}$, as defined in \cite[Section~6.2]{La87}, because they, together with~$j$, generate the field $\Q(\zeta_N)(X(N))$. As explained in \cite[bottom of page~66]{La87}, the  $q$-expansion of $f_{r,s}$ is of the form
\begin{equation}
\label{eqexpweb}
\sum_{k=k_0}^\infty a_{k,r}(\zeta_N^s)q^{k/N},
\end{equation}
where  ${a_{k,r}(t)\in \Q[t]}$ are polynomials depending on~$k$ and~$r$, but not on~$s$. 

For ${\sigma=(\begin{smallmatrix}1&0\\0&d\end{smallmatrix})}$ we have ${f_{r,s}^\sigma=f_{r,sd}}$. Hence the $q$-expansion of $f_{r,s}^\sigma$ is like~\eqref{eqexpweb}, but with $\zeta_N^s$ replaced by $\zeta_N^{sd}$. Using item~\ref{idetsigma} we find
$$
a_{k,r}(\zeta_N^{sd})= a_{k,r}((\zeta_N^\sigma)^s)= a_{k,r}(\zeta_N^s)^\sigma,
$$
as wanted. 
\end{proof}

Let~$G$ be a subgroup of $\GL_2(\Z/N\Z)$ containing ${-I}$. We denote by~$X_G$ the associated modular curve. It corresponds to the $G$-invariant subfield of the field ${\Q(\zeta_N)\bigl(X(N)\bigr)}$.  The constant subfield of this field is ${\Q(\zeta_N)^{\det G}}$, where ${\det:\GL_2(\Z/N\Z)\to(\Z/N\Z)^\times}$ is the determinant,  and we identify $(\Z/N\Z)^\times$ with the Galois group $\gal(\Q(\zeta_N)/\Q)$. In particular, if  ${\det G=(\Z/N\Z)^\times}$ then the constant subfield is~$\Q$ and  the corresponding modular curve~$X_G$ is defined (that is, has a geometrically irreducible model) over~$\Q$.

For a subgroup~$H$ of $(\Z/N\Z)^\times$ put
\begin{equation}
\label{egh}
G_H=\{g\in G: \det g\in H\}. 
\end{equation}
In particular, ${G_{(\Z/N\Z)^\times}=G}$ and ${G_1=G\cap\SL_2(\Z/N\Z)}$. If~$H$ is contained in $\det G$, then the subfield of $\Q(\zeta_N)\bigl(X(N)\bigr)$ stabilized by~$G_H$ is $K(X_G)$, where ${K=\Q(\zeta_N)^H}$. 

\begin{remark}
\label{rorbits}
Let~$M_N$ be the subset of the abelian group $(\Z/N\Z)^2$ consisting of the elements of exact order~$N$. Then the set of cusps of the modular curve $X_G$ is in natural one-to-one correspondence with the set 
${G_1\backslash M_N}$ 
of orbits of the natural left action of~$G_1$ on $M_N$ \cite[Lemma~2.3]{BI11}. 

The cusps are defined over the cyclotomic field $\Q(\zeta_N)$. 
Identifying  the groups ${\gal\bigl(\Q(\zeta_N)/\Q\bigr)}$ and  ${(\Z/N\Z)^\times}$, the natural left action of ${(\Z/N\Z)^\times}$ on the set ${G_1\backslash M_N}$ coincides with the Galois action on the cusps. Hence, if~$H$ is a subgroup of $(\Z/N\Z)^\times$ then the set of $H$-orbits of cusps stands in a one-to-one correspondence with left $G_H$-orbits on $M_N$. 
\end{remark}

\subsection{The nearest cusp} 
\label{ssnear}

Let~$\Gamma$ be the subgroup of ${\Gamma(1)=\SL_2(\Z)}$ obtained by lifting~$G_1$. Then the set of complex points $X_G(\C)$ is analytically isomorphic to ${\Gamma\backslash\bar\HH}$, where ${\bar\HH=\HH\cup\Q\cup\{i\infty\}}$ is the extended Poincaré plane. Similarly, $Y_G(\C)$ is  analytically isomorphic to ${\Gamma\backslash\HH}$. 

We denote by~$\calF$ the familiar fundamental domain of the modular group ${\Gamma(1)=\SL_2(\Z)}$: the open hyperbolic triangle with vertices ${e^{i\pi/3},e^{2i\pi/3},i\infty}$, together with the geodesics ${[i,e^{i\pi/3}]}$  and ${[e^{i\pi/3},i\infty)}$. There is a natural bijection ${Y(1)(\C)\leftrightarrow \calF}$, and the image of ${P\in Y(1)(\C)}$ under this bijection will be denoted ${\tau(P)}$.
More generally, the image of ${P\in Y_G(\C)}$ under the map ${Y_G(\C)\to Y(1)(\C)\to\calF}$ will also be denoted $\tau(P)$. Alternatively, we can define $\tau(P)$ as the single ${\tau\in\calF}$ with the property 
\begin{equation}
\label{ejptau}
j(P)=j(\tau),
\end{equation}
where we use the convention for~$j$ from Section~\ref{ssconv}.

We also consider the slightly smaller set
$$
\ocalF=\{z\in \calF: |z|>1\}. 
$$
In other words,~$\ocalF$ is~$\calF$ with the geodesic ${[i,e^{i\pi/3}]}$ removed.

For every ${\sigma \in \Gamma(1)}$ we define  the set ${\calF(\sigma)\subset Y_G(\C)}$   as the image of $\sigma\calF$ in ${Y_G(\C)=\Gamma\backslash\HH}$; similarly, $\ocalF(\sigma)$ is the  image of $\sigma\ocalF$. Clearly 
$$
\ocalF(\sigma)=\{P\in \calF(\sigma): |\tau(P)|>1\}. 
$$

The sets $\calF(\sigma)$ and $\ocalF(\sigma)$ depend only on the coset $\Gamma\sigma$; in particular, there are exactly ${[\Gamma(1):\Gamma]}$ distinct sets $\calF(\sigma)$. They are pairwise disjoint and cover $Y_G(\C)$:
\begin{align*}
\bigcup_{\Gamma\sigma} \calF(\sigma)=Y_G(\C), \qquad \calF(\sigma)\cap \calF(\sigma')=\varnothing \quad (\Gamma\sigma\ne \Gamma\sigma'),
\end{align*}
the union being over the cosets ${\Gamma\backslash\Gamma(1)}$.

Next, for every cusp~$c$ we define ${\Omega_c\subseteq X_G(\C)}$ and ${\oOmega_c\subset \Omega_c}$ by
\begin{equation}
\label{eomegac}
\Omega_c = \bigcup_{\sigma(i\infty)=c} \calF(\sigma)\cup\{c\}, \qquad \oOmega_c = \bigcup_{\sigma(i\infty)=c} \ocalF(\sigma)\cup\{c\}=\{P\in \Omega_c: |\tau(P)|>1\},
\end{equation}
the union being over all ${\sigma	 \in \Gamma}$ such that $\sigma(i\infty)$ represents the cusp~$c$. 

This can be made more explicit as follows. Let ${e=e_c}$ be the ramification index of the branch cover ${X_G\to X(1)}$  at~$c$. Fix some ${\sigma	 \in \Gamma}$ such that $\sigma(i\infty)$ represents the cusp~$c$, and define, for ${k\in \Z}$,
\begin{equation}
\label{esigmak}
\sigma_k=\sigma\circ \begin{pmatrix}
1&k\\0&1
\end{pmatrix} .
\end{equation}
Then 
\begin{equation}
\label{odecomp}
\Omega_c=\bigcup_{k=0}^{e-1}\calF(\sigma_k)\cup \{c\}, \qquad 
\oOmega_c=\bigcup_{k=0}^{e-1}\ocalF(\sigma_k)\cup \{c\}.
\end{equation}

The sets~$\Omega_c$
 are pairwise disjoint and  cover $X_G(\C)$:
$$
\bigcup_{c}\Omega_c= X_G(\C), \qquad \Omega_c\cap\Omega_{c'}=\varnothing \quad (c\ne c'). 
$$
If ${P\in X_G(\C)}$ belongs to~$\Omega_c$, we call~$c$ the \textsl{nearest cusp} to~$P$. Note that the set $\oOmega_c$ is open in the complex topology, and is the maximal open subset of $\Omega_c$. 

In practical calculations, we select a full set~$\Sigma$ of representatives of cosets ${\Gamma\backslash\Gamma(1)}$. Then we have 
$$
X_G(\C)=\bigcup_{\sigma\in \Sigma}\calF(\sigma). 
$$
We develop methods for finding integral points on each $\ocalF(\sigma)$. Our methods do not work for the points~$P$ with ${|\tau(P)|=1}$, but these can be found just by checking all integral values of~$j$ from~$0$ to~$1728$, see Proposition~\ref{preal} below. 

We build~$\Sigma$ as follows: first, for every cusp~$c$ we pick ${\sigma_c\in \Gamma(1)}$ such that ${\sigma_c(i\infty)}$ represents~$c$. After this is done, we define~$\Sigma$ as the set of all $\sigma_{c,k}$ with ${0\le k\le e_c-1}$ for every~$c$, where ${\sigma_{c,k}=\sigma_c\circ \bigl(\begin{smallmatrix}
1&k\\0&1
\end{smallmatrix}\bigr)}$. See Section~\ref{sscusps} for a concrete example.

\subsection{The $q$-parameter at a cusp}
\label{ssqpar}

For ${P\in \Omega_c}$ we define the $q$-parameter $q_c(P)$ by
${q_c(P)=e^{2\pi i\tau(P)}}$, with the convention ${\tau(c)=i\infty}$ and ${q_c(c)=0}$. 
This $q_c$ is a holomorphic function on  $\oOmega_c$. We have
\begin{equation}
\label{ejpq}
j(P)=j(q_c(P))
\end{equation}
(see the same convention on~$j$ as above). 
Since ${\im\,\tau(P)\ge\sqrt3/2}$ %for ${P\in \Omega_c}$, 
we have  
\begin{equation}
\label{eless}
|q_c(P)|\le e^{-\pi\sqrt3}<0.0044 \qquad (P\in \Omega_c). 
\end{equation}

As in Section~\ref{ssnear} denote by ${e=e_c}$ the ramification index at~$c$ of  ${X_G\to X(1)}$. Then $q_c^{1/e_c}$ can be viewed as a ``local analytic parameter'' at~$c$. This means the following: if ${u\in \C(X_G)}$ is a $\C$-rational function on~$X_G$, then in a neighborhood of~$c$ we have
\begin{equation*}
%\label{elocalp}
\log|u(P)|= \frac{\ord_cu}{e_c}\log|q_c(P)|+O(1).
\end{equation*}

This can also be expressed in terms of Taylor expansions. Loosely speaking, it means that for ${P\in \oOmega_c}$ and for a suitable choice of the $e$th root $q_c(P)^{1/e}$ we have
\begin{equation}
\label{etaylorloose}
u(P) = \varsigma_cq_c(P)^{\ord_cu/e}+O\bigl(|q_c(P)|^{(\ord_cu+1)/e}\bigr),
\end{equation}
where $\varsigma_c$ (which is well-defined up to multiplication by an $e$th root of unity) and the implied constant depend only on~$u$, but not on~$P$. We make this more precise as follows. 

Fix ${\sigma\in \Gamma(1)}$ such that ${\sigma(i\infty)}$ represents the cusp~$c$, and define~$\sigma_k$ as in~\eqref{esigmak}.     There exists a non-zero complex number ${\varsigma_c=\varsigma_{c,\sigma}}$  
such that the following holds. For   ${P\in \oOmega_c}$ define ${q_c(P)^{1/e}= e^{2\pi i\tau(P)/e}}$. Then 
\begin{equation}\label{etaylorexact}
u(P)= \varsigma_c\bigl(e^{2k\pi i/e}q_c(P)^{1/e}\bigr)^{\ord_c(u)}+ O\bigl(|q_c(P)|^{(\ord_cu+1)/e}\bigr) \qquad (P \in \ocalF(\sigma_k), \ k\in \Z). 
\end{equation}
Due to  decomposition~\eqref{eomegac}, this gives an exact version of the ``Taylor expansion''~\eqref{etaylorloose} on $\oOmega_c$. 

Note that ${\varsigma_c=\varsigma_{c,\sigma}}$  does depend on the choice of~$\sigma$; if we change~$\sigma$ then~$\varsigma_c$ would be multiplied be an $e$th root of unity.

\subsection{More on~$j$}

The following property will be routinely used.

\begin{proposition}
\label{preal}
For a non-cusp point ${P\in X_G(\C)}$ the following two conditions are equivalent.

\begin{enumerate}
\item
${j(P)\in \R}$;

\item
${\re(\tau(P))\in \{0,1/2\}}$ or ${|\tau(P)|=1}$. 
\end{enumerate}
More precisely:
\begin{equation*}
\begin{array}{lclcl}
j(P) \in [1728,+\infty) &\Longleftrightarrow& \re(\tau(P))=0&\Longleftrightarrow&q_c(P)>0,\\
j(P) \in (-\infty,0] &\Longleftrightarrow& \re(\tau(P))=1/2&\Longleftrightarrow& q_c(P) <0,\\
j(P) \in [0,1728] &\Longleftrightarrow& |\tau(P)|=1,
\end{array}
\end{equation*}
where~$c$ is the nearest cusp to~$P$. 
\end{proposition}

\begin{proof}
It suffices to show that for ${\tau\in \calF}$ we have 
\begin{alignat}3
\label{etwelvecube}
j(\tau)&\in [1728,+\infty)\quad  &&\Longleftrightarrow&& \quad \re(\tau)=0,\\
\label{enegzero}
j(\tau) &\in (-\infty,0]\quad  &&\Longleftrightarrow&& \quad \re(\tau)=1/2,\\
\label{emodone}
j(\tau) &\in [0,1728]\quad  &&\Longleftrightarrow&& \quad |\tau|=1.
\end{alignat}
%Since~$j$ takes on~$\calF$ every complex value exactly once, it suffices to prove only the ``$\Leftarrow$'' implications. 
Since the coefficients in the $q$-expansion of~$j$ are real, we have 
${j(iy)\in \R}$ for ${y>0}$.  Using 
$$
j(i)=1728, \qquad  \lim_{y\to+\infty}j(iy)=+\infty,
$$
this proves~\eqref{etwelvecube}, because~$j$ takes on~$\calF$ every complex value exactly once.  In a similar fashion one proves~\eqref{enegzero}: we have ${j(1/2+iy)\in \R}$ and 
$$
j(e^{\pi i/3})=0, \qquad \lim_{y\to+\infty}j(1/2+iy)=-\infty. 
$$
Finally, using again that the coefficients in the $q$-expansion of~$j$ are real, we deduce that, for ${\tau \in \HH}$, we have
${j(-\bar\tau)=\overline{j(\tau)}}$. 
Together with ${j(-1/\tau)=j(\tau)}$ this implies that ${j(\tau)\in \R}$ when ${|\tau|=1}$. Since 
${j(e^{\pi i/3})=0}$ and ${j(i)=1728}$,
this proves~\eqref{emodone}. 
\end{proof}

We shall also need an approximate formula for the $j$-invariant. Write the  $q$-expansion as\footnote{The coefficients~$c_n$ cannot be confused with the cusps.}
$$
j(q)=c_{-1}q^{-1}+c_0+c_1q+c_2q^2+\ldots, 
$$
with ${c_{-1}=1}$,  ${c_0=744}$, ${c_1=196884}$ etc. For a non-negative integer~$N$ write
$$
j_N(q)=  \sum_{n=-1}^{N-1} c_nq^n.
$$
In particular, ${j_0(q)=q^{-1}}$. 

\begin{lemma}
For ${P\in \Omega_c}$ we have
\begin{equation}
\label{eestj}
j(P) = j_N(q_c(P))  +R_N, \qquad |R_N|\le j(e^{-\pi\sqrt3})-j_N(e^{-\pi\sqrt3}) 
\end{equation}
for any non-negative integer~$N$. In particular, for ${N=0}$ we have 
\begin{equation}
\label{e2079}
\bigl|j(P) - q_c(P)^{-1}\bigr|\le 2079.
\end{equation}
\end{lemma} 

\begin{proof}
Since~$j$ is $\Gamma(1)$-invariant, we may assume that~$c$ is the cusp at infinity and ${q_c(P)=q(P)}$. Since the coefficients~$c_n$ are known to be positive and ${|q(P)|\le e^{-\pi\sqrt3}}$, we have
$$
|j(P) - j_N(q_c(P))|\le \sum_{n=N}^\infty c_n|q(P)|^n \le \sum_{n=N}^\infty c_n\bigl|e^{-\pi\sqrt3}\bigr|^n= j(e^{-\pi\sqrt3})-j_N(e^{-\pi\sqrt3}),
$$
proving~\eqref{eestj}. 
In particular, for ${N=0}$ we obtain 
$$
\bigl|j(P) - q_c(P)^{-1}\bigr|\le j(e^{-\pi\sqrt3}) - e^{\pi\sqrt3} < 2309.6- 230.7 <  2079. 
$$
\end{proof}

Note that the exact value of $j(e^{-\pi\sqrt3})$ is available: $j(e^{-\pi\sqrt3})= 40500  (35010 -20213\sqrt{3})$. 

\section{Integral points and Baker's method on modular curves}
\label{sbakergen}

In this section we give a general overview of Baker's method applied to modular curves. For more details, see~\cite{BS13,Bi02} and Sha's thesis~\cite{Sh13}. 

Let~$N$ and~$G$ be as in Section~\ref{sxg}, let~$K$ be a number field containing ${\Q(\zeta_N)^{\det G}}$ and $\OO_K$  the ring of integers of~$K$. %(This notation cannot be confused with the notation~$\OO$ for the $G$-orbits in the previous section.) 
We define the set of integral points 
$$
X_G(\OO_K)=\{P\in X_G(K): j(P)\in \OO_K\}.
$$
Recall that the the \textsl{height} of ${\alpha\in \OO_K}$ is defined by 
$$
\height(\alpha) = [K:\Q]^{-1}\sum_{\sigma:K\hookrightarrow \C}\log^+|\alpha^\sigma|, \qquad \log^+=\max\{\log,0\}. 
$$
the sum being over the complex embeddings of~$K$.

%One can similarly define the set of $S$-integral points for a finite set~$S$ of places of~$K$. 
%Denote by $\height(\alpha)$ the usual absolute logarithmic height of an algebraic number~$\alpha$. 
We want to bound the height  $\height(j(P))$ for ${P\in X_G(\OO_K)}$. We show how to do this under the assumption 
\begin{equation}
\label{e3}
\nu_\infty(G)\ge 3,
\end{equation}
where $\nu_\infty(G)$ denotes the number of cusps of $X_G$. 

A  \textsl{modular unit} is  a rational function (defined over~$\bar K$) on~$X_G$ with no zeros and no poles outside the cusps. Equivalently, ${u\in \bar K(X_G)}$ is a modular unit if both~$u$ and~$u^{-1}$ are integral over the ring ${\Q[j]}$. Principal divisors of modular units form a subgroup in the group of degree~$0$ divisors supported on the cusps. The latter is a free abelian group of rank ${\nu_\infty(G)-1}$, so the group of principal divisors of modular units must be of rank not exceeding ${\nu_\infty(G)-1}$. It is of fundamental importance for us that it is of the maximal possible rank; this is sometimes called the \textsl{Manin--Drinfeld theorem}.

\begin{theorem}%[Manin-Drinfeld]
\label{tmandr}
The principal divisors of modular units form a free abelian group of rank ${\nu_\infty(G)-1}$. 
\end{theorem}
See \cite[Chapter~4, Theorem~2.1]{La95}.
Here is an immediate consequence.

\begin{corollary}
\label{cmandr}
Assume that ${\nu_\infty(G)\ge 3}$. Then for any cusp~$c$ there exists a non-constant modular unit~$u$ such that ${u(c)=1}$. 
\end{corollary}

If ${j(P)\in \OO_K}$ then 
${\height(j(P))= [K:\Q]^{-1}\sum_{\sigma:K\hookrightarrow \C}\log^+|j(P)^\sigma|}$,
the sum being over the complex embeddings of~$K$. 
For some embedding~$\sigma$ we have 
${\height(j(P)) \le \log|j(P)^\sigma|}$.
We fix this embedding from now on and view~$K$ as a subfield of~$\C$. Thus, we have to bound $|j(P)|$ from above.  

The point~$P$ belongs to one of the sets~$\Omega_c$, defined in~\eqref{eomegac}, and the corresponding~$c$ is the ``nearest cusp'' to~$P$. Now, since ${\nu_\infty(G)\ge 3}$, we may use Corollary~\ref{cmandr} and find a non-constant modular unit~$u$ with ${u(c)=1}$. The rational function~$u$ is defined over the number field $K(\zeta_N)$. 

If ${u(P)=1}$ then it is easy to bound~$P$ as one of the zeros of the rational function ${u-1}$. From now on we assume that ${u(P)\ne 1}$.  Since ${u(c)=1}$, we have 
$$
u(P)= 1+ O(|q_c(P)|^{1/{e_c}}).
$$
(Here and below in this section, the constant implied by the $O(\cdot)$-notation, as well as by the Vinogradov notation ``$\ll$'' and ``$\gg$'', may depend on~$N$ and~$K$, but not on~$P$.) Thus, $u(P)$ is a complex algebraic number, distinct from~$1$ but ``close'' to~$1$ if $q_c(P)$ is small.

Since both~$u$ and $u^{-1}$ are integral over $\Q[j]$,  there exist non-zero  ${A_1, A_2\in \Z}$, which can be easily determined explicitly, such that $A_1u$ and $A_2u^{-1}$ are integral over $\Z[j]$. Since ${j(P)\in \OO_K}$, both ${A_1u(P)}$ and $A_2u(P)^{-1}$ belong to $\OO_{K(\zeta_N)}$. It follows that there are only finitely many possibilities for the principal ideal $(u(P))$ (viewed as a fractional ideal in the field $K(\zeta_N)$). 

Fixing a system $\eta_1, \ldots, \eta_r$ of fundamental units of $K(\zeta_N)$, we obtain 
${u(P)=\eta_0\eta_1^{b_1}\cdots \eta_r^{b_r}}$, where $\eta_0$ belongs to a finite subset of~$K$ (that can be explicitly determined), and  ${b_1,\ldots, b_r}$ are rational integers depending on~$P$. We obtain the inequality
\begin{equation}
\label{eforb}
\left|\eta_0\eta_1^{b_1}\cdots \eta_r^{b_r} -1\right|\ll q_c(P)^{1/{e_c}}.
\end{equation}
Let $B=\max(|b_1|,\ldots,|b_r|)$.
It is easy to show that ${B\ll \height(\eta)}$, see \cite[bottom of page~77]{Bi02}. It follows that ${B\ll \height(u(P))+1}$. On the other hand, the general property of quasi-equivalence of heights on an algebraic curve implies that ${\height(u(P)) \ll \height(j(P))+1}$. It follows that 
\begin{equation}
\label{equasi}
B\ll \height(j(P)) \le \log|j(P)| = \log|q_c(P)^{-1}| +O(1). 
\end{equation}

On the other hand, one can bound the  left-hand side of~\eqref{eforb} from below using the so-called \textsl{Baker's inequality}, which implies that either the  left-hand side of~\eqref{eforb} is~$0$ (in which case ${u(P)=1}$ and  $\height(j(P))$ is bounded), or it is bounded from below by ${\exp(-\kappa\log \max(B,3))}$, where~$\kappa$ is a positive effective constant depending on ${\eta_0, \eta_1, \ldots, \eta_r}$ but independent of~$B$. Combining this with~\eqref{eforb}, we obtain the estimate
${\log|q_c(P)^{-1}| \ll \log \max(B,3)}$. Together with~\eqref{equasi} this bounds $|q_c(P)|$ away from zero, which  implies a bound for $|j(P)|$ from above. See~\cite{BS13}, where this approach is used to bound explicitly integral points on $X_\ns^+(p)$ for ${p\ge 7}$.

In a similar fashion one can study $S$-integral points on $X_G$: the new ingredients to be added  are the $p$-adic version of Baker's inequality, due to Yu~\cite{Yu07}, and the $p$-adic analogue of the notion of the ``nearest cusp'', see \cite[Section~3]{BP11}. We do not go into this in the present article.

To make the argument above explicit, one needs to construct modular units explicitly. The standard tool for this is \textsl{Siegel functions}, see Section~\ref{ssiegel} below.  One also needs explicit version for various statements above like the quasi-equivalence of heights, etc. All this can be found in the Ph.D. thesis of Sha~\cite{Sh13,Sh14}.

In the present work, we are interested in a somewhat different task: not just to bound the heights of integral points, but to determine them completely. We restrict ourselves to the case ${K=\Q}$ and ${N=p}$ a prime number. In this case the most interesting class of modular curves for which integral points are unknown is  $X_\ns^+(p)$, when the group~$G$ is the normalizer of a non-split  Cartan subgroup of $\GL_2(\F_p)$. 

The principal point here is that bounding the height of integral points, even explicitly in all parameters, is not sufficient for the actual calculation of the points. The problem is that the bounds obtained by Baker's method are excessively huge and not suitable for direct enumeration. 

Fortunately, one can reduce Baker's bound using the technique of numerical Diophantine approximation introduced by Baker and Davenport~\cite{BD69}. This reduction is described in detail in \cite{BH96,BH99,Ha00} in the context of the Diophantine equation of Thue. 
Recall that this is the equation of the form 
${f(x,y)=A}$,
where the ${f(x,y)\in \Z[x,y]}$  is a $\Q$-irreducible form of degree ${n\ge 3}$, and~$A$  is a non-zero integer. In~\cite{BH98} the method was extended to the superelliptic Diophantine equations.  Here we adapt this reduction method to  modular curves.

Several observations are to be made.

\begin{enumerate}

\item
Usually, to perform the computations, one  should know explicitly the algebraic data of the  number field(s) involved (in the case of Thue equation, this is the  field generated over~$\Q$ by a root of $f(1,y)$). By the algebraic data we mean here the unit group (with explicit generators), the class group (again, for every class one should have an explicit ideal representing this class), and so on. Fortunately, in the special case of the curve $X_\ns^+(p)$ these tasks are radically simplified.

First, the field we are going to deal with is the real cyclotomic field ${\Q(\zeta_p+\bar\zeta_p)}$  (or a subfield, see below) for which the unit group (or at least a full-rank subgroup of the latter, which is sufficient, see below) is given explicitly by the \textsl{circular units}. To be precise, the index of the group of circular units in the full unit group of ${\Q(\zeta_p+\bar\zeta_p)}$ is equal to the real class number $h_p^+$; see, for instance, \cite[Theorem~8.2]{Wa97}.

In the range ${p<100}$ that we are working, we have ${h_p^+=1}$, see the recent article of Miller~\cite{Mi15}, who extended the earlier work of Masley~\cite{Mas78}. Hence in this range circular units form  the full unit group.

Second, the only ideal we are going to deal with is the one above~$p$, which is principal and has an obvious explicit generator ${(\zeta_p-\bar\zeta_p)^2}$. This was already used in~\cite{BHV01} for solving Thue equations ${\Phi_n(x,y)=p}$, where $\Phi_n(1,y)$ is the $n$-th real cyclotomic polynomial, and~$p$ is a prime divisor of~$n$. 

\item
To make the calculations more efficient, it is in some cases useful to replace the field ${\Q(\zeta_p+\bar\zeta_p)}$ by a smaller subfield, if possible.  %Hence, instead of using the group~$G_1$, as in Proposition~\ref{pmdr}, one should use~$G_H$ with a suitably chosen subgroup ${H\ni-1}$ of $(\Z/p\Z)^\times$. In this case we would be able to replace ${\Q(\zeta_p+\bar\zeta_p)}$ by the smaller field ${\Q(\zeta_p)^H}$. 
This was suggested in~\cite{BH99} and was very efficiently exploited in~\cite{BHV01}.

In the setting of the present paper, with this trick we reduced the running for $p=97$ by a factor of~$4.5$ using the subfield of degree $16$.
We also tried to use the degree $12$ subfield, however then the sieves become too imprecise so that it becomes more expensive to exclude potential lattice points on the curve~$\gamma$ (cf. Section~\ref{secEllipsoidsSieve}), and we ended up with a running time improvement by only a factor of~$1.5$.
%And for the degree $24$ subfield, computing the fundamental units did not finish within two days and we decided to interrupt their computation.

\item
In principle, it is not necessary to have the full unit group; a full-rank subgroup would suffice, as explained in~\cite{Ha00}. In particular, since circular units form a full rank subgroup for all~$p$, our method must work for all primes, not only for those where the circular units form the full unit group. Note that 
this was used already in~\cite{BHV01} (in a different context).

In the present work we use only full unit groups as they were always computable for our parameters, either due to the above-mentioned result of Miller~\cite{Mi15}, or using the \textsf{PARI} package. 
However, one should keep this opportunity in mind for further applications.

\item
Adapting numerical methods developed for Thue equations to modular
curves is not straightforward. In the Thue case one  has formulas with
very strong error estimates, typically $O(|x|^{-n})$, where~$n$ is the
degree of the equation; see, for instance
\cite[Proposition~2.4.1]{BH96}.  This is quite good even for small
solutions~$x$.

However, for modular curves of level~$p$ we have, typically, errors
$O(|j(P)|^{-1/p})$.
Larger errors mean that we have to check more false candidates to find
all solutions.
Therefore in Section~\ref{secEllipsoidsSieve} we considerably improve the involved sieves.
If the error bounds are too weak, one can use higher order asymptotic
expansions for the modular functions involved, see
Appendix~\ref{appendix}.
At the point where this becomes computationally too expensive, we stop
the sieve and start an extra search, which checks all $j$'s with small
modulus separately.
\end{enumerate}

\section{Siegel functions}
\label{ssiegel}
In this section we  recall the principal facts about Klein forms and Siegel functions.  For more details the reader can consult \cite[Section~2.1]{KL81} and~\cite{KS10}. We call a positive integer~$N$  a \textit{denominator} of ${a\in\Q}$ if ${Na\in\Z}$. For instance, $2020$ is a denominator of $1/2$.

\subsection{Klein forms and Siegel functions}
\label{ssklesie}
Let ${\tilbfa=(\tila_1,\tila_2)\in \Q^2}$ be such that ${\tilbfa\notin \Z^2}$. We denote by $\gerk_\tilbfa(\tau)$ the Klein form associated to~$\tilbfa$, which is a holomorphic function on the Poincaré plane~$\HH$. We collect some properties of Klein forms in the proposition below.

\begin{proposition}
\label{pprokl}
\begin{enumerate}

\item
The Klein forms do not vanish on~$\HH$. 

\item
\label{iklega1}
The Klein forms behave well under the action of $\Gamma(1)$: for  ${\sigma=(\begin{smallmatrix}a&b\\c&d\end{smallmatrix})\in \Gamma(1)}$   we have
\begin{equation*}
%\label{eklega1}
\gerk_\tilbfa\circ\sigma(\tau)=(c\tau+d)^{-1}\gerk_{\tilbfa\sigma}(\tau),
\end{equation*}
where  ${\sigma(\tau)= \frac{a\tau+b}{c\tau+d}}$. In particular, with ${\sigma=-I}$ this gives
\begin{equation}
\label{eklemin}
\gerk_{-\tilbfa}=-\gerk_\tilbfa.
\end{equation}

\item
\label{iklemod1}
For ${\tilbfa=(\tila_1,\tila_2) \in \Q^2\setminus\Z^2}$ and ${\tilbfb=(\tilb_1,\tilb_2)\in \Z^2}$ we have 
\begin{equation*}
\gerk_{\tilbfa+\tilbfb}= \eps(\tilbfa,\tilbfb)\gerk_\tilbfa, \qquad \eps(\tilbfa,\tilbfb)= (-1)^{\tilb_1\tilb_2+\tilb_1+\tilb_2}e^{\pi i(\tila_1\tilb_2-\tila_2\tilb_1)}.
\end{equation*}
Notice that ${\eps(\tilbfa,\tilbfb)^{2N}=1}$, where~$N$ is a denominator of~$\tilbfa$ (a common denominator of~$\tila_1$ and~$\tila_2$).

\item
Let~$N$ be a denominator of~$\tilbfa$.   Then~$\gerk_\tilbfa$ is  ``nearly'' $\Gamma(N)$-automorphic of weight~$-1$. Precisely, for
${\sigma=(\begin{smallmatrix}a&b\\c&d\end{smallmatrix}) \in \Gamma(N)}$ we have 
$$
\gerk_\tilbfa\circ\sigma(\tau)=\eps'(\tilbfa,\sigma)(c\tau+d)^{-1}\gerk_\tilbfa(\tau), \qquad  \eps'(\tilbfa,\sigma)^{2N}=1. 
$$

\end{enumerate}

\end{proposition}

The following result is a consequence of the properties above. 

\begin{proposition}
\label{pkle}
Let~$N$ be a denominator of~$\tilbfa$. Then $\gerk_\tilbfa^{2N}$ depends only on the residue class of~$\tilbfa$ modulo $\Z^2$, and it is $\Gamma(N)$-automorphic of weight ${-2N}$. 
\end{proposition}

Further, for ${\tilbfa=(\tila_1,\tila_2)\in \Q^2\smallsetminus \Z^2}$ we define the \textit{Siegel function}  $g_\tilbfa(\tau)$ by
$$
g_\tilbfa(\tau)= \gerk_\tilbfa(\tau)\eta(\tau)^2,
$$ 
where $\eta(\tau)$ is the Dedekind $\eta$-function. 

Since ${\eta(\tau)^{24}=\Delta(\tau)}$ is $\Gamma(1)$-automorphic of weight~$12$, Proposition~\ref{pkle} implies the following.

\begin{theorem}
\label{t12N}
In the set-up of Proposition~\ref{pkle}, the function $g_\tilbfa^{12N}$ depends only on the residue class of~$\tilbfa$ modulo $\Z^2$, and is $\Gamma(N)$-automorphic of weight~$0$. 
\end{theorem}

It follows, in particular,  that  Siegel functions~$g_\tilbfa$  are algebraic over the field $\C(j)$ (because so are $\Gamma(N)$-automorphic functions).   In addition to this,~$g_\tilbfa$ 
is holomorphic and does not vanish on the Poincaré plane~$\HH$ (because so are the Klein forms and the Dedekind~$\eta$). It follows that both~$g_\tilbfa$ and $g_\tilbfa^{-1}$  must be integral over the ring $\C[j]$. Actually, a stronger 
assertion holds (see, for instance, Proposition~2.2 from~\cite{BP11a}).

\begin{proposition}
\label{psiu}
Let~$N$ be the smallest denominator of~$\tilbfa$  and~$\zeta_N$ a primitive $N$-th root of unity. Then both~$g_\tilbfa$ and   ${\left(1-\zeta_N\right)g_\tilbfa^{-1}}$ are integral over  $\Z[j]$. 
\end{proposition}

\subsection{An approximate formula}

As usual, write ${q=q(\tau)=e^{2\pi i\tau}}$. 
For a rational number~$a$ we define ${q^a=e^{2\pi i a\tau}}$. Then the Siegel function~$g_\tilbfa$ has the following infinite product presentation \cite[page~29]{KL81}:
\begin{equation}
\label{epga}
g_\tilbfa(\tau)= -q^{B_2(\tila_1)/2}e^{\pi i\tila_2(\tila_1-1)}\prod_{n=0}^\infty(
1-q^{n+\tila_1}e^{2\pi i\tila_2})\left(1-q^{n+1-\tila_1}e^{-2\pi i \tila_2}\right), 
\end{equation}
where ${B_2(T)=T^2-T+{1}/{6}}$ is the second Bernoulli polynomial. Together with Proposition~\ref{pprokl}:\ref{iklemod1}, this implies that 
\begin{equation}
\label{eordqga_and_ella}
\ord_qg_\tilbfa = \ell_\tilbfa := B_2(\tila_1-\lfloor \tila_1\rfloor)/2.
\end{equation} 
Here the $q$-order $\ord_q$ is defined by ${\lim_{q\to0}q^{\ord_qg_\tilbfa}g_\tilbfa(q)\ne0,\infty}$. 

In fact, we have the following quantitative statement. 

\begin{proposition}
\label{pellaass}
Put
\begin{equation}
\label{egammaa}
\varrho_\tilbfa=
\begin{cases}
-e^{\pi i\tila_2(\tila_1-1)}, & \tila_1\ne 0,\\
-e^{\pi i\tila_2(\tila_1-1)}(1-e^{2\pi i \tila_2}), &\tila_1=0.
\end{cases}
\end{equation}
Then, for a non-zero ${\tilbfa=(\tila_1,\tila_2) \in \Q^2\cap[0,1)^2}$ and ${\tau\in \HH}$ we have
\begin{equation}
\label{ehardestimate}
\left|\log\frac{g_\tilbfa(\tau)}{\varrho_\tilbfa q^{\ell_\tilbfa}}\right| \le 
\begin{cases}
\frac1{1-|q|}\left(\frac{|q|^{\tila_1}}{1-|q|^{\tila_1}}+\frac{|q|^{1-\tila_1}}{1-|q|^{1-\tila_1}}\right),& \tila_1>0,\\
\frac{2|q|}{(1-|q|)^2}.& \tila_1=0.
\end{cases}
\end{equation}
\end{proposition}

\begin{proof}
For ${|z|<1}$  we have 
${\bigl|\log (1+z)\bigr|  \le |z|/(1-|z|)}$. 
Hence, when ${\tila_1>0}$, we can bound the terms of the product expansion~\eqref{epga}  as
\[
\bigl|\log(1-q^{n+\tila_1}e^{2\pi i\tila_2})\bigr|\le   \frac{|q|^{n+\tila_1}}{1-|q|^{\tila_1}}, \qquad 
\bigl|\log(1-q^{n+1-\tila_1}e^{-2\pi i\tila_2})\bigr|\le\frac{|q|^{n+1-\tila_1}}{1-|q|^{1-\tila_1}}. 
\]
Adding this up for all ${n\geq 0}$, we obtain
\[
\left|\log\frac{g_\tilbfa(\tau)}{\varrho_\tilbfa q^{\ell_\tilbfa}}\right| \le \sum_{n=0}^\infty \left(\frac{|q|^{n+\tila_1}}{1-|q|^{\tila_1}} +\frac{|q|^{n+1-\tila_1}}{1-|q|^{1-\tila_1}}\right)
=\frac1{1-|q|}\left(\frac{|q|^{\tila_1}}{1-|q|^{\tila_1}}+\frac{|q|^{1-\tila_1}}{1-|q|^{1-\tila_1}}\right),
\]
which proves~\eqref{ehardestimate} in the case ${\tila_1>0}$. 

In the case ${\tila_1=0}$ we re-write~\eqref{epga} as
$$
g_\tilbfa(\tau)= \varrho_\tilbfa q^{\ell_\tilbfa}\prod_{n=1}^\infty\left(
1-q^{n}e^{2\pi i\tila_2}\right)\left(1-q^{n}e^{-2\pi i \tila_2}\right).
$$
We bound 
$$
\bigl|\log(1-q^{n}e^{2\pi i\tila_2})\bigr|, \bigl|\log(1-q^{n}e^{-2\pi i\tila_2})\bigr|\le \frac{|q|^n}{1-|q|} \qquad (n\ge 1)
$$
Adding this up for ${n\ge 1}$, we prove~\eqref{ehardestimate} in the case ${\tila_1=0}$ as well. 
\end{proof}

\begin{corollary}
\label{cellaass}
In the set-up of Proposition~\ref{pellaass} assume that ${\tau \in \calF}$. Let~$N$ be  a denominator of~$\tila_1$, and assume that ${N\ge 5}$. Then   
\begin{equation}
\label{eeasyestimate}
\left|\log\frac{g_\tilbfa(\tau)}{\varrho_\tilbfa q^{\ell_\tilbfa}}\right| \le N
|q|^{1/N}. 
\end{equation}
\end{corollary}

\begin{proof}
Assume first that ${\tila_1>0}$, in which case  ${\tila_1, 1-\tila_1\ge1/N}$. The function 
$$
x\mapsto \frac{x}{1-x}+\frac{|q|x^{-1}}{1-|q|x^{-1}}
$$
is increasing on the interval ${[|q|^{1/2}, 1)}$. Hence 
$$
\frac{|q|^{\tila_1}}{1-|q|^{\tila_1}}+\frac{|q|^{1-\tila_1}}{1-|q|^{1-\tila_1}}\le \frac{|q|^{1/N}}{1-|q|^{1/N}}+\frac{|q|^{1-1/N}}{1-|q|^{1-1/N}}
$$
Since ${\tau \in \calF}$, we have ${|q|\ge e^{-\pi\sqrt3}}$. Using the assumption ${N\ge 5}$, we obtain 
\begin{align*}
 &\frac{|q|^{1/N}}{1-|q|^{1/N}}\le \frac{|q|^{1/N}}{1-e^{-\pi\sqrt3/N}}=\frac{e^{\pi\sqrt3/N}}{e^{\pi\sqrt3/N}-1}|q|^{1/N}\le {e^{\pi\sqrt3/5}}\frac N{\pi\sqrt3}|q|^{1/N},\\
 &\frac{|q|^{1-1/N}}{1-|q|^{1-1/N}}\le \frac{|q|^{3/5}}{1-|q|^{4/5}}|q|^{1/N}\le \frac{e^{-3\pi\sqrt3/5}}{1-e^{-4\pi\sqrt3/5}}|q|^{1/N}\le \frac N5\frac{e^{-3\pi\sqrt3/5}}{1-e^{-3\pi\sqrt3/5}}|q|^{1/N}.
\end{align*}
It follows that  
$$
\left|\log\frac{g_\tilbfa(\tau)}{\varrho_\tilbfa q^{\ell_\tilbfa}}\right| \le \frac{1}{1-e^{-\pi\sqrt3}}\left( \frac {e^{\pi\sqrt3/5}}{\pi\sqrt3}+\frac{e^{-3\pi\sqrt3/5}}{5(1-e^{-4\pi\sqrt3/5})}\right)N|q|^{1/N} <0.56N|q|^{1/N}, 
$$
which proves~\eqref{eeasyestimate}, in a stronger form, in the case ${\tila_1>0}$. 

The case ${\tila_1=0}$ is much simpler: 
$$
\left|\log\frac{g_\tilbfa(\tau)}{\varrho_\tilbfa q^{\ell_\tilbfa}}\right| \le \frac{2|q|}{(1-|q|)^2}\le \frac{2e^{-4\pi\sqrt3/5}}{(1-e^{-\pi\sqrt3})^2}|q|^{1/N} \le 
0.03|q|^{1/N}, 
$$
which is much better than~\eqref{eeasyestimate}.
\end{proof}

More refined approximate formulas can be found in the Appendix, see Proposition~\ref{papprsiegel} therein.

\subsection{Simplest modular units}
\label{sssmu}
Now let us fix a positive integer~$N$. We have the natural group isomorphism ${(N^{-1}\Z/\Z)^2\cong (\Z/N\Z)^2}$, and, with some abuse of speech, we identify the two groups. In particular, for ${\bfa\in (\Z/N\Z)^2}$ we have the corresponding element in ${(N^{-1}\Z/\Z)^2}$,  and for this latter we may fix a lifting  ${\tilbfa\in N^{-1}\Z^2}$, which will be called a lifting of~$\bfa$ to ${N^{-1}\Z^2}$.

By Theorem~\ref{t12N}, for ${\bfa\in N^{-1}\Z^2\smallsetminus\Z^2}$ 
the function $$u_\bfa := g_\tilbfa^{12N}$$  does not depend on a particular choice of the lifting~$\tilbfa$ and defines a $\C$-rational function on the 
modular curve $X(N)$. 
Identity~\eqref{eklemin} implies that $u_\bfa =u_{-\bfa}$. 

The infinite product~\eqref{epga}  implies that the $q$-expansion of~$u_\bfa$ has coefficients in the cyclotomic field $\Q(\zeta_N)$. 
By~\cite[Proposition 6.9(1)]{Shi71}, it follows that ${u_\bfa\in \Q(\zeta_N)\bigl(X(N)\bigr)}$. 
Moreover, the Galois action of the group $\GL_2(\Z/N\Z)$ on the field ${\Q(\zeta_N)\bigl(X(N)\bigr)}$ (see Section~\ref{sxg}) coincides with the action induced by the natural right action of $\GL_2(\Z/N\Z)$ on the set ${(\Z/N\Z)^2}$   in the following sense: for a non-zero  ${\bfa \in (\Z/N\Z)^2}$  and ${\sigma \in \GL_2(\Z/N\Z)}$ we have 
\begin{equation}
\label{egalua}
u_{\bfa\sigma}=u_\bfa^\sigma.
\end{equation}
See \cite[Section~4.2]{BP11} for more details.

The functions~$u_\bfa$ give the simplest explicit examples of the \textit{modular units}, already mentioned in Section~\ref{sbakergen}: they have no zeros and no poles outside the cusps. It follows that their principal divisors generate a free abelian subgroup of rank at most ${\nu_\infty(N)-1}$, where $\nu_\infty(N)$ is the number of cusps of $X(N)$. It turns out that this rank is the maximal possible, which provides an explicit form of the Manin--Drinfeld theorem (Theorem~\ref{tmandr}):

\begin{theorem}
\label{tmdr}
The principal divisors $(u_\bfa)$ generate a free abelian group of rank ${\nu_\infty(N)-1}$. 
\end{theorem}
For the proof see Theorem~3.1 in 
\cite[Chapter~2]{KL81}.

In fact, one can show that already the principal divisors $(u_\bfa)$, where~$\bfa$ runs through the set~$M_N$, consisting of the elements of $(\Z/N\Z)^2$ of exact order~$N$, generate a free abelian group of rank ${\nu_\infty(N)-1}$. The number of such~$\bfa$ is $2\nu_\infty(N)$. It follows that, besides the relations ${u_\bfa=u_{-\bfa}}$, there can exist exactly one relation between the principal divisors $(u_\bfa)$ with ${\bfa\in M_N}$. This relation is 
$$
\sum_{\bfa\in M_N}(u_\bfa)=0.
$$ 
In fact, we have a more precise statement.

\begin{lemma}
\label{lprodofus}
In the above set-up we have 
\begin{equation}
\label{eprodua}
\prod_{\bfa \in M_N}u_\bfa= \pm\Phi_N(1)^{12N},
\end{equation}
where~$\Phi_N(t)$ is the $N$-th cyclotomic polynomial. In particular, if ${N=p}$ is a prime number, we have
\begin{equation}
\label{eproduap}
\prod_{\bfa \in M_p}u_\bfa= \pm p^{12p}.
\end{equation}
\end{lemma}
%(One can show that the sign is actually~$+$.)
This will be used in the proof of the principal relation (Section~\ref{sprinrel}) that our algorithm is based upon.

\begin{proof}
%Let us prove~\eqref{eprodua}. 
Since the set~$M_N$ is stable with respect to ${\GL_2(\Z/N\Z)}$, the left-hand side of~\eqref{eprodua} is stable with respect to the Galois action over the field $\Q(X(1))$. Hence it is a unit on the curve $X(1)$, defined over~$\Q$. 
Since $X(1)$ has only one cusp, it has no non-constant units. Hence the left-hand side of~\eqref{eprodua} is a constant belonging to~$\Q$.

To determine the value of this constant, we evaluate it at the cusp at infinity. For each ${\bfa \in (\Z/N\Z)}$ we choose the lifting ${\tilbfa=(\tila_1,\tila_2)\in \Q^2}$ such that ${0\le a_1,a_2<1}$. The left-hand side of~\eqref{eprodua} is a product of a root of unity and the terms of the type ${\bigl(1-e^{2\pi i a_2}q^{n+\tila_1}\bigr)^{12N}}$ and of the type ${\bigl(1-e^{2\pi i -a_2}q^{n+1-\tila_1}\bigr)^{12N}}$, where~$n$ runs through non-negative integers, and $(\tila_1,\tila_2)$ runs through the liftings of the elements of the set $M_N$.   When we set ${q=0}$, all these terms become~$1$ except the terms ${\bigl(1-e^{2\pi i a_2}q^{n+a_1}\bigr)^{12N}}$ with ${n=0}$ and ${a_1=0}$. Hence, up to a root of unity,  the left-hand side of~\eqref{eprodua} is
$$
\prod_{\genfrac{}{}{0pt}{}{a_2\in N^{-1}\Z/\Z}{\text{$a_2$ is of order~$N$}}} \bigl(1-e^{2\pi i a_2}\bigr)^{12N}= \prod_{\genfrac{}{}{0pt}{}{0\le k <N}{(k,N)=1}}\bigl(1-e^{2\pi i k/N}\bigr)^{12N}= \Phi_N(1)^{12N}.
$$
Since the only roots of unity in~$\Q$ are~$\pm1$, this proves~\eqref{eprodua} and the lemma.  
\end{proof}

\section{General modular units}
\label{squad}

In this section we review and complement some of the results of Kubert and Lang~\cite{KL81}. 
Our purpose is to construct ``economical'' modular units on the curve $X_G$. 

The ``naive''  approach is as follows. Let~$G$ be a subgroup of ${\GL_2(\Z/N\Z)}$ and~$H$ a subgroup of $\det G$, which  itself is a subgroup in  $(\Z/N\Z)^\times$, viewed as the Galois group of the cyclotomic field $\Q(\zeta_N)$.  Then~$H$ left-acts naturally on the set of the cusps of~$X_G$. Denote by $\nu_\infty(G)$ the number of cusps and by $\nu_\infty(G, H)$  the number of  $H$-orbits of cusps.

On the other hand, the group~$G_H$, defined in~\eqref{egh}, right-acts on the set $(\Z/N\Z)^2$.  If ${\OO\subset (\Z/N\Z)^2}$ is a non-zero orbit of this action, then 
\begin{equation}
\label{euo}
\prod_{\bfa\in \OO}u_\bfa 
\end{equation}
is a rational function on the curve $X_G$ defined over the field ${\Q(\zeta_N)^H}$. 

It is not difficult to deduce from Theorem~\ref{tmdr} that 
the principal divisors defined by products~\eqref{euo}, where~$\OO$ runs the non-zero $G_H$-orbits,  generate a free abelian group whose rank is  ${\nu_\infty(G, H)-1}$. 

Product~\eqref{euo} can be written as 
\begin{equation}
\label{ega12N}
\prod_{\bfa\in \OO}g_\tilbfa^{12N},
\end{equation}
where ${\tilbfa \in N^{-1}\Z^2}$ is a lifting of ${\bfa\in (\Z/N\Z)^2}$, as defined in Section~\ref{sssmu}. The principal goal of this section is constructing rational functions on~$X_G$ of the form 
\begin{equation*}
%\label{ega12N}
\prod_{\bfa\in \OO}g_\tilbfa^{m},
\end{equation*}
where~$m$ is much smaller than $12N$, which   is crucial for numerical purposes. %This is the principal goal of this section. 

\subsection{Quadratic relations}
\label{ssquad}

Let~$N$ be a positive integer. As in Section~\ref{sssmu},
we identify the groups ${(N^{-1}\Z/\Z)^2}$ and~${(\Z/N\Z)^2}$, which allows us to lift every ${\bfa\in (\Z/N\Z)^2}$ to some  ${\tilbfa \in N^{-1}\Z^2}$. By a lifting of a set ${A\subset (\Z/N\Z)^2}$ we mean a mapping ${A\to N^{-1}\Z^2}$ such that for every ${\bfa \in A}$ its image ${\tilbfa\in N^{-1}\Z^2}$ is a lifting of~$\bfa$ in the sense defined above.

Our principal tool will be the following result 
of Kubert and Lang~\cite{KL81}, see Theorem~5.2 in Chapter~3. 

\begin{theorem}
\label{tquad}

To every non-zero ${\bfa=(a_1,a_2) \in (\Z/N\Z)^2}$ we associate an integer $m(\bfa)$.  Fix a lifting ${\bfa\mapsto\tilbfa}$ of the set of non-zero elements of ${(\Z/N\Z)^2}$. % ${\tilbfa\in N^{-1}\Z^2}$ of~$\bfa$.  %again denoted by~$\bfa$. 
Put 
\begin{equation}
\label{esumma}
\Lambda=\sum_{\genfrac{}{}{0pt}{}{\bfa \in (\Z/N\Z)^2}{\bfa \ne 0}}m(\bfa).
\end{equation}
   
\begin{enumerate}
\item
Assume that~$N$ is odd. Then 
\begin{equation}
\label{eprodkama}
\prod_{\genfrac{}{}{0pt}{}{\bfa \in (\Z/N\Z)^2}{\bfa \ne 0}}\gerk_\tilbfa^{m(\bfa)} 
\end{equation}
is $\Gamma(N)$-automorphic (of weight $-\Lambda$) if and only if 
\begin{equation}
\label{equadN}
\sum_{\genfrac{}{}{0pt}{}{\bfa \in (\Z/N\Z)^2}{\bfa \ne 0}}m(\bfa)a_1^2= 
\sum_{\genfrac{}{}{0pt}{}{\bfa \in (\Z/N\Z)^2}{\bfa \ne 0}}m(\bfa)a_2^2= 
\sum_{\genfrac{}{}{0pt}{}{\bfa \in (\Z/N\Z)^2}{\bfa \ne 0}}m(\bfa)a_1a_2= 0. 
\end{equation}
\item
\label{iprodgama}
Assume that ${\gcd(N,6)=1}$. 
Then the function 
\begin{equation}
\label{eprodgama}
\prod_{\genfrac{}{}{0pt}{}{\bfa \in (\Z/N\Z)^2}{\bfa \ne 0}}g_\tilbfa^{m(\bfa)} 
\end{equation}
is $\Gamma(N)$-automorphic (of weight~$0$) if and only if~\eqref{equadN} holds and ${12\mid \Lambda}$.  
\end{enumerate}
\end{theorem}

\begin{remark}
\label{rquad}
\begin{enumerate}
\item
Kubert and Lang call~\eqref{equadN} ``quadratic relations'' (modulo~$N$). 
\item
One may note that 
\begin{equation}
\label{ewithdelta}
\prod_{\genfrac{}{}{0pt}{}{\bfa \in (\Z/N\Z)^2}{\bfa \ne 0}}g_\tilbfa^{m(\bfa)}=\prod_{\genfrac{}{}{0pt}{}{\bfa \in (\Z/N\Z)^2}{\bfa \ne 0}}\gerk_\tilbfa^{m(\bfa)} \cdot \Delta^{\Lambda/12},
\end{equation}
where ${\Delta=\eta^{24}}$. 
\item
The assumption ${\gcd(N,6)=1}$ is purely technical: in a slightly modified form the statement holds 
true when~$N$ is divisible by~$2$ and/or by~$3$. However, assuming that ${\gcd(N,6)=1}$ will not hurt us, since we shall apply Theorem~\ref{tquad} only when~$N$ is prime and ${N\ge 7}$. 

\item
\label{izetan}
Theorem~\ref{tquad} implies that product~\eqref{eprodgama} defines a function ${u\in \C\bigl(X(N)\bigr)}$. By considering the $q$-expansion, as in Section~\ref{sssmu}, we conclude that in fact ${u\in \Q(\zeta_N)\bigl(X(N)\bigr)}$. 
\end{enumerate}
\end{remark}

Contrary to product~\eqref{ega12N}, product~\eqref{eprodgama} may depend on the choice of the lifting ${\bfa\mapsto\tilbfa}$. Proposition~\ref{pprokl}:\ref{iklemod1} implies that if we choose a different lifting ${\bfa\mapsto\tilbfa'}$ then~\eqref{eprodkama} and~\eqref{eprodgama} will be multiplied by a $2N$-th root of unity. Though this is pretty trivial, we state this as a proposition for further reference.

\begin{proposition}
\label{ptriv}
For every non-zero ${\bfa \in (\Z/N\Z)^2}$ pick an integer $m(\bfa)$ and fix  \emph{two} liftings ${\bfa\mapsto\tilbfa}$ and ${\bfa\mapsto \tilbfa'}$ of the set of non-zero elements of ${(\Z/N\Z)^2}$. Then there exists a $2N$-th root of unity~$\eps$ such that 
\begin{equation}
\label{etrivkle}
\prod_{\genfrac{}{}{0pt}{}{\bfa \in (\Z/N\Z)^2}{\bfa \ne 0}}\gerk_{\tilbfa'}^{m(\bfa)} =\eps\prod_{\genfrac{}{}{0pt}{}{\bfa \in (\Z/N\Z)^2}{\bfa \ne 0}}\gerk_\tilbfa^{m(\bfa)}.
\end{equation}
If additionally ${12\mid \Lambda}$, where~$\Lambda$ defined in~\eqref{esumma}, then 
\begin{equation}
\label{etrivg}
\prod_{\genfrac{}{}{0pt}{}{\bfa \in (\Z/N\Z)^2}{\bfa \ne 0}}g_{\tilbfa'}^{m(\bfa)} =\eps\prod_{\genfrac{}{}{0pt}{}{\bfa \in (\Z/N\Z)^2}{\bfa \ne 0}}g_\tilbfa^{m(\bfa)}.
\end{equation}
If  $2\mid m(\bfa)$ for every $\bfa$ then 
\begin{equation}
\label{etriveps}
\eps^N=1.
\end{equation}
\end{proposition}

\begin{proof}
Statements~\eqref{etrivkle} and~\eqref{etriveps} follow from Proposition~\ref{pprokl}:\ref{iklemod1}, and~\eqref{etrivg} follows from~\eqref{etrivkle} and~\eqref{ewithdelta}.
\end{proof}

\subsection{Galois action}
As we mentioned in Section~\ref{sssmu}, the Galois action by the group $\GL_2(\Z/N\Z)$ on the ``simplest'' modular units ${u_\bfa=g_\tilbfa^{12N}}$ is very easy to describe: it is given by relation~\eqref{egalua}. We want to obtain a similar result for ``general'' modular units~\eqref{eprodgama}.

\begin{proposition}
\label{psiegal}
Assume the set-up of Theorem~\ref{tquad}:\ref{iprodgama}, so that 
$$
u=\prod_{\genfrac{}{}{0pt}{}{\bfa \in (\Z/N\Z)^2}{\bfa \ne 0}}g_\tilbfa^{m(\bfa)} 
$$
defines a function in $\Q(\zeta_N)(X_G)$ (see Remark~\ref{rquad}:\ref{izetan}). 

\begin{enumerate}
\item
\label{isl2}
Assume that ${\sigma \in \SL_2(\Z/N\Z)}$ and let~$\tilsigma$ be a lifting of~$\sigma$ to $\Gamma(1)$. Then
\begin{equation}
\label{efsigma}
u^\sigma=\prod_{\genfrac{}{}{0pt}{}{\bfa \in (\Z/N\Z)^2}{\bfa \ne 0}}g_{\tilbfa\tilsigma}^{m(\bfa)}. 
\end{equation}

\item
\label{igl2}
Assume that ${\sigma \in \GL_2(\Z/N\Z)}$. Then it has a lifting ${\tilsigma\in \mathrm M_2(\Z)}$ such that~\eqref{efsigma} holds.

\end{enumerate}
\end{proposition}

\begin{proof}
The first part is a consequence of Proposition~\ref{pgalomod}:\ref{ifsig}, Proposition~\ref{pprokl}:\ref{iklega1} and~\eqref{ewithdelta}. Indeed, write ${\tilsigma=(\begin{smallmatrix}a&b\\c&d\end{smallmatrix})}$, and recall that ${\Delta=\eta^{24}}$ is $\Gamma(1)$-automorphic of weight~$12$. We obtain
\begin{align*}
u^\sigma(\tau)&=u\circ\tilsigma(\tau)\\
&= \prod_{\genfrac{}{}{0pt}{}{\bfa \in (\Z/N\Z)^2}{\bfa \ne 0}}\bigl(\gerk_\tilbfa\circ\tilsigma(\tau)\bigr)^{m(\bfa)} \cdot \bigl(\Delta\circ\tilsigma(\tau)\bigr)^{\Lambda/12}\\
&=(c\tau+d)^{-\Lambda}\prod_{\genfrac{}{}{0pt}{}{\bfa \in (\Z/N\Z)^2}{\bfa \ne 0}}\gerk_{\tilbfa\tilsigma}(\tau)^{m(\bfa)}\cdot  (c\tau+d)^{\Lambda} \Delta(\tau)^{\Lambda/12}\\
&=\prod_{\genfrac{}{}{0pt}{}{\bfa \in (\Z/N\Z)^2}{\bfa \ne 0}}g_{\tilbfa\tilsigma}(\tau)^{m(\bfa)},
\end{align*}
as wanted.

In the proof of the second part, we may assume that~$\sigma$ is of the form $(\begin{smallmatrix}1&0\\0&d\end{smallmatrix})$, because  any ${\sigma\in \GL_2(\Z/N\Z)}$ can be presented as $\sigma_1\sigma_2$ with ${\sigma_1 \in \SL_2(\Z/N\Z)}$ and~$\sigma_2$ of this form. We lift ${\sigma=(\begin{smallmatrix}1&0\\0&d\end{smallmatrix})}$ as ${\tilsigma=(\begin{smallmatrix}1&0\\0&\tild\end{smallmatrix})}$, and the result follows immediately from Proposition~\ref{pgalomod}:\ref{igalqexp} and the infinite product~\eqref{epga}.
\end{proof}

\subsection{Economical modular units on $X_G$}
\label{ssecon}
In this section to avoid technicalities we restrict to  prime level. Thus, let ${p\ge 5}$ be a prime number,~$G$  a subgroup in $\GL_2(\F_p)$ and~$H$  a subgroup in $\det G$. 
The group~$G_H$, defined in~\eqref{egh}, right-acts on the set
${M_p=\F_p^2\smallsetminus \{0\}}$ (as in the previous section, we tacitly identify the sets $\F_p^2$ and ${(p^{-1}\Z/\Z)^2}$). 
Let ${\OO\subset M_p}$  be an orbit of this action, or, more generally, a $G_H$-invariant subset of~$M_p$.  We fix a lifting  ${\bfa\mapsto\tilbfa}$ of the set~$M_p$ (as defined in the beginning of Section~\ref{ssquad})  and we want to find an exponent~$m$ such that 
\begin{equation}
\label{euoo}
u=\prod_{\bfa\in \OO}g_\tilbfa^m
\end{equation}
defines a function in $K(X_G)$, where ${K=\Q(\zeta_p)^H}$. Clearly,
${m=12p}$ would do. It turns out that in some cases one can do much better, sometimes introducing a root of unity factor. We fix a $p$-th primitive root of unity and denote it by $\zeta_p$. 

\begin{theorem}
\label{tecoun}
Let ${p\ge 5}$ be a prime number and ${G\ni-I}$  a  
subgroup of $\GL_2(\F_p)$ such that that $|G|$ is not divisible by~$p$. 
Let~$H$ be a subgroup of $\det G$ and ${\OO\subset M_p}$  a $G_H$-invariant subset of~$M_p$ satisfying
\begin{equation}
\label{equadoo}
\sum_{\bfa\in \OO}a_1^2= \sum_{\bfa\in \OO}a_1a_2=\sum_{\bfa\in \OO}a_2^2=0. 
\end{equation}
Let~$m$ be an integer such that 
\begin{equation}
\label{e12moo}
2\mid m, \qquad 12\mid m|\OO|.
\end{equation} 
Fix a lifting  ${\bfa\mapsto\tilbfa}$ of the set~$\OO$   and define~$u$ as in~\eqref{euoo}.  Then~$u$ defines a function in $\Q(\zeta_p)(X_G)$ (denoted by~$u$ as well). Further, 
%when the liftings~$\tilbfa$ are chosen in a suitable way (explained in the proof) we have  
there exists ${k\in \Z}$ (which is unique $\bmod\, p$ when ${H\ne 1}$) such that  
${\zeta_p^ku\in K(X_G)}$, where ${K=\Q(\zeta_p)^H}$. 
\end{theorem}

The proof requires a  lemma, which is the simplest special case of the Kummer theory (see any textbook in algebra). 

\begin{lemma}
\label{lkum}
Let~$p$ be a prime number and~$F$ a field of characteristic distinct from~$p$. Let~$\alpha$ be an element in the algebraic closure~$\bar F$, and ${\zeta_p\in \bar F}$ a primitive $p$-th root of unity. Assume that ${\alpha^p\in F}$. Then either ${[F(\alpha):F]=p}$ or there exists  ${k\in \Z}$ (which is unique $\bmod\, p$ when ${\zeta_p\notin F}$) such that ${\zeta_p^k\alpha \in F}$. In particular, if ${\zeta_p\in F}$ then either ${[F(\alpha):F]=p}$ or ${\alpha\in F}$. 
\end{lemma}

\begin{proof}[Proof of Theorem~\ref{tecoun}]
Theorem~\ref{tquad} together with Remark~\ref{rquad}:\ref{izetan} imply that~$u$ defines a function in $\Q(\zeta_p)\bigl(X(p)\bigr)$. We want to study the Galois action of $G_H$ on~$u$. %We claim that ${f^p \in K(X_G)}$. To see this, 
Thus, fix ${\sigma \in G_H}$. 
Proposition~\ref{psiegal}:\ref{igl2} implies that there exists a lifting ${\tilsigma\in \mathrm M_2(\Z)}$ such that 
\begin{equation}
\label{efsigoo}
u^\sigma= \prod_{\bfa\in \OO} g_{\tilbfa\tilsigma}^m. 
\end{equation}
Since~$\OO$ is $G_H$-invariant, we have ${\OO\sigma^{-1}=\OO}$. Consider a different lifting ${\bfa\mapsto\tilbfa'}$ of ~$\OO$ defined by  ${\tilbfa'= \widetilde{\bfa\sigma^{-1}}\tilsigma}$, where ${\widetilde{\bfa\sigma^{-1}}}$ is the  lifting of ${\bfa\sigma^{-1}}$. Then~\eqref{efsigoo} can be rewritten as 
$$
u^\sigma= \prod_{\bfa\in \OO} g_{\tilbfa'}^m.
$$
Now Proposition~\ref{ptriv} implies that ${u^\sigma/u}$ is a $p$-th root of unity. We have proved that~$u^p$ is invariant under the Galois action by~$G_H$, which implies that ${u^p \in K(X_G)}$, the $G_H$-invariant subfield of $\Q(\zeta_p)\bigl(X(p)\bigr)$. Since the degree
$
{\bigl[\Q(\zeta_p)\bigl(X(p)\bigr):K(X_G)\bigr]=|G_H|}
$
is not divisible by~$p$ (because $|G|$ is not divisible by~$p$ by the assumption),
Lemma~\ref{lkum} completes the proof.
\end{proof}

There is an important special case when~$u$ itself belongs to $K(X_G)$, without multiplication by a root of unity. Assume that~$G_H$ contains $(\begin{smallmatrix}1&0\\0&-1\end{smallmatrix})$. In this case ${\bfa=(a_1,a_2)}$ belongs to a $G_H$-orbit~$\OO$ if and only if its ``complex conjugate'' ${\bar\bfa=(a_1,-a_2)}$ does. We say that a lifting ${\bfa\mapsto\tilbfa}$ \textsl{respects complex conjugation} if the following holds: if ${\bfa=(a_1,a_2)\in \OO}$ is lifted to ${\tilbfa=(\tila_1,\tila_2)}$, then the lifting of~$\bar\bfa$ is $(\tila_1,-\tila_2)$. This can be expressed briefly as ${\tilde{\bar\bfa}=\bar\tilbfa}$. 

\begin{corollary}
\label{cecoun}
In the set-up of Theorem~\ref{tecoun} assume that ${(\begin{smallmatrix}1&0\\0&-1\end{smallmatrix})\in G_H}$ and that the lifting respects complex conjugation. Then ${u\in K(X_G)}$.
\end{corollary}

\begin{proof}
The assumption ${(\begin{smallmatrix}1&0\\0&-1\end{smallmatrix})\in G_H}$ implies that ${K\subseteq \Q(\zeta_p+\bar\zeta_p)}$. Further, since the lifting respects complex conjugation, we have ${u^\iota=u}$, where ${\iota=(\begin{smallmatrix}1&0\\0&-1\end{smallmatrix})}$. The subfield of $\Q(\zeta_p)(X_G)$ stabilized by~$\iota$ is ${\Q(\zeta_p+\bar\zeta_p)(X_G)}$. Thus, ${u\in \Q(\zeta_p+\bar\zeta_p)(X_G)}$ and ${\zeta_p^ku\in K(X_G)}$ with ${K\subseteq \Q(\zeta_p+\bar\zeta_p)}$. It follows that ${\zeta_p^k}$ lies in  ${\Q(\zeta_p+\bar\zeta_p)}$ which is only possible if ${\zeta_p^k=1}$.
\end{proof}

\subsection{An approximate formula}

Using Proposition~\ref{pellaass} and Corollary~\ref{cellaass}, we may obtain   approximate expressions for the modular units constructed in Section~\ref{ssecon}. Let~$p$,~$G$,~$H$ and~$\OO$ be as in Section~\ref{ssecon}. In particular, as in Theorem~\ref{tecoun}, we will assume that 
\begin{equation*}
%\label{egndivp}
\text{$|G|$ is not divisible by~$p$}.
\end{equation*}
 Let~$c$ be a cusp of $X_G$. We define the sets~$\Omega_c$,~$\oOmega_c$ and the $q$-parameter $q_c$ as in Sections~\ref{ssnear} and~\ref{ssqpar}, that is, ${q_c(P) = e^{2\pi i\tau(P) }}$ for ${P\in \Omega_c}$. We also fix ${\sigma\in \Gamma(1)}$ such that ${\sigma(i\infty)}$ represents the cusp~$c$ and define~$\sigma_k$ as in~\eqref{esigmak}.

Since the ramification of ${X(p)\to X(1)}$ at all cusps is~$p$, the ramification of ${X_G\to X(1)}$ at~$c$  is either~$1$ or~$p$. Moreover, since $|G|$ is not divisible by~$p$,  the ramification is~$p$ at all cusps. 

As in the previous section, we fix a lifting ${\bfa\mapsto\tilbfa}$ of the set~$M_p$.

For a subset ${A\in M_p}$   define the quantities
\begin{equation}
\label{elara}
\ell_A=\sum_{\bfa\in A}\ell_\tilbfa, \qquad \varrho_A=\prod_{\bfa\in A}\varrho_\tilbfa, 
\end{equation}
where~$\ell_\tilbfa$ is defined in~\eqref{eordqga_and_ella} and~$\varrho_\tilbfa$ is defined in~\eqref{egammaa}. 
Note that~$\ell_A$ is independent of the fixed lifting, but~$\varrho_A$ depends on it and is well-defined only up to multiplication by a $p$th root of unity. However, we will mainly deal with the absolute value $|\varrho_A|$,  
which is independent of the lifting. 

\begin{proposition}
\label{passu}
Assume that~$G$ is a semi-simple subgroup of $\GL_2(\F_p)$  and that~$\OO$ and~$m$ satisfy the hypothesis of Theorem~\ref{tecoun}.  
Define~$u$ as  in~\eqref{euoo}, and define ${\varsigma_c=\varsigma_{c,\sigma} }$ as in Section~\ref{ssqpar}. Then 
\begin{equation}
\label{eordcoef}
\frac{\ord_cu}{p}=m\ell_{\OO\sigma}, \qquad \varsigma_c= \varrho_{\OO\sigma}^m\cdot(\text{a $p$th root of unity}) .
\end{equation}
Furthermore, for ${P\in \calF(\sigma_k)}$ we have
\begin{equation}
\label{esimpleestimateulog}
\log|u(P)|= m\ell_{\OO\sigma} \log|q_c(P)|+m\log|\varrho_{\OO\sigma}|+O_1\bigl( 
pm|\OO||q_c(P)|^{1/p}\bigr). 
\end{equation}
\end{proposition} 

\begin{proof}
Replacing~$G$ by ${\sigma^{-1}G\sigma}$ and~$\OO$ by $\OO\sigma$, we may assume that~$c$ is represented by $i\infty$. We may further assume, without loss of generality, that ${\sigma=I}$. In this special case Proposition~\ref{passu} follows immediately from Proposition~\ref{pellaass} and Corollary~\ref{cellaass} applied with ${N=p}$ to every~$\tilbfa$ lifting some ${\bfa\in \OO}$; recall that ${p\ge 5}$ by the assumption.  
\end{proof}

In the sieving algorithm of Section~\ref{secEllipsoidsSieve} we use  more refined approximate formulas from Appendix~\ref{appendix}.

\subsection{An example}
We conclude this section with an example. It will not be used in the sequel, but it gives a good illustration of how Theorem~\ref{tecoun} can be used. 

We take as~$G$ the diagonal subgroup of $\GL_2(\F_p)$ and set ${H=\{1,-1\}}$, so that
$$
G_H=\bigl\{(\begin{smallmatrix}a&0\\0&d\end{smallmatrix}): ad=\pm1\bigr\} 
$$
and ${K=\Q(\zeta_p+\bar\zeta_p)}$.

The right $G_H$-action on~$M_p$ has   ${(p-1)/2}$ distinct orbits. They are of the form ${\{\bfa: a_1a_2=\pm c\}}$ with ${c=1, \ldots, (p-1)/2}$. The quadratic relations~\eqref{equadoo} are clearly satisfied, and to have~\eqref{e12moo} it suffices to take 
\begin{equation*}
m=
\begin{cases}
2, & p\equiv 1\bmod3,\\
6, & p\equiv -1\bmod3.
\end{cases}
\end{equation*} 
Selecting a lifting respecting the complex conjugation, we obtain ${(p-1)/2}$ modular units in the field $K(X_G)$.

\section{Cusp points and units on  \texorpdfstring{$X_\ns^+(p)$}{Xns+(p)}}
\label{scpu}

From now on we restrict to the case when ${N=p}$ is a prime number and~$G$ is  the normalizer of a non-split Cartan subgroup of $\GL_2(\Z/p\Z)$. A very detailed account of various properties of this curve  (even for an arbitrary~$N$) can be found in Sections~3 and~6 of Baran's article~\cite{Ba10}. 

We may and will assume that
\begin{equation}
\label{eexpl}
G=\left\lbrace \begin{pmatrix} \alpha & \Xi\beta \\ \beta & \alpha \end{pmatrix} , \begin{pmatrix} \alpha & \Xi\beta \\ -\beta & -\alpha \end{pmatrix} : \alpha, \beta \in \F_p,\ (\alpha,\beta)\ne (0,0) \right\rbrace ,
\end{equation}
where~$\Xi$ is a quadratic non-residue modulo~$p$, which will be fixed from now on. In particular, one can take ${\Xi=-1}$  if ${p\equiv 3\bmod4}$. 

We fix until the end of the article a lifting ${\bfa\mapsto\tilbfa}$ of the set $M_p$ to $p^{-1}\Z^2$, which respects complex conjugation (as defined before Corollary~\ref{cecoun}) and which has in addition to this the following property:
\begin{equation}
\label{e01}
\text{if ${\tilbfa=(\tila_1,\tila_2)}$ is a lifting of ${\bfa\in M_p}$ then ${0\le \tila_1<1}$.}
\end{equation}

\subsection{Cusps}
\label{sscusps}
The curve ${X_G=X_\ns^+(p)}$ has ${(p-1)/2}$ cusps, defined over the real cyclotomic fields ${\Q(\zeta_p+\bar\zeta_p)}$, and the Galois group ${\gal\bigl(\Q(\zeta_p+\bar\zeta_p)/\Q\bigr)=\F_p^\times/\{\pm1\}}$ acts transitively on the cusps. 

According to Remark~\ref{rorbits}, the  cusps stay in one-to-one correspondence with
the  orbits of the left $G_1$-action on the set ${M_p=\F_p^2\smallsetminus\{(0,0)\}}$. These orbits are the sets defined by ${x^2-\Xi y^2=\pm c}$, where~$c$ runs through representatives of cosets  ${\F_p^\times/\{\pm1\}}$, the cusp at infinity corresponding to ${c=1}$. 

For every ${c\in\F_p^\times/\{\pm1\}}$ fix ${(a,b) \in \F_p^2}$ such that ${a^2-\Xi b^2=c^{-1}}$ and let $\sigma_c$ be a lifting of the matrix $\left(\begin{smallmatrix}ca&b\Xi\\cb&a\end{smallmatrix}\right)$ to~$\Gamma(1)$. For ${c=1}$ we take ${(a,b)=(1,0)}$ and ${\sigma_1=I}$. Then the set ${\bigl\{\sigma_c(i\infty): c\in\F_p^\times/\{\pm1\}\bigr\}}$ is a full system of representatives of cusps on~$\bar\HH$, and the set
$$
\Sigma=\left\{\sigma_c\circ\begin{pmatrix}
1&k\\0&1
\end{pmatrix}:c\in\F_p^\times/\{\pm1\}, \ k=0,\ldots, p-1\right\}
$$
is a complete  system of representatives of cosets of $\Gamma_\ns^+\backslash\Gamma(1)$. This is a special case of the construction  explained in Section~\ref{ssnear}.  

In the sequel we fix a subgroup~$H$ of $\F_p^\times$ containing~${-1}$ 
and put ${d=[\F_p^\times:H]}$. 
In particular,
$$
d= [K:\Q],
$$
where ${K=\Q(\zeta_p)^H}$. The group~$H$ acts on the set of cusps by Galois conjugation, and this action has exactly~$d$ orbits, each of them being defined over~$K$ as a set. The Galois group ${\gal(K/\Q)=\F_p^\times/H}$ acts on the set of $H$-orbits transitively.  These $H$-orbits of cusps are in one-to-one correspondence with the sets defined by ${x^2-\Xi y^2\in cH}$, with $cH$ running through the cosets $\F_p^\times/H$.

\subsection{Units}
\label{ssuoo}

Besides the left action, the group~$G_H$ acts on the set~$M_p$ from the right. There are again~$d$  orbits of this action, and they are defined by ${\Xi x^2- y^2\in cH}$.    These orbits will be used to define modular units in $K(X_G)$. Recall that we fixed a lifting ${\bfa\mapsto\tilbfa}$  of~$M_p$ to ${p^{-1}\Z^2}$, respecting the complex conjugation.

\begin{theorem}
\label{tu}
Let~$\OO$ be a right $G_H$-orbit on~$M_p$. Pick a lifting ${\bfa \mapsto\tilbfa}$ of~$\OO$ to ${p^{-1}\Z^2}$. 
Put
\begin{equation}
\label{edefm}
m=
\begin{cases}
2, & \text{if $3\mid (p+1)|H|$}, \\
6, & \text{otherwise}. 
\end{cases}
\end{equation} 
Then the product
\begin{equation}
\label{efum}
u_\OO=\prod_{\bfa\in \OO}g_\tilbfa^m
\end{equation}
is well-defined (it depends only on the orbit~$\OO$ but not on the particular lifting) and it defines a function in $K(X_G)$. 
\end{theorem}

We deduce this theorem from Theorem~\ref{tecoun} (more precisely from Corollary~\ref{cecoun}) using  some elementary lemmas about finite fields. We thank 
 Julia Baoulina for useful explanations and for the proof of Lemma~\ref{lba} below. 
 
\begin{lemma}
\label{lln}
Let  ${P(x_1, \ldots, x_n)\in \F[x_1, \ldots, x_n]}$ be a polynomial over a finite field ${\F=\F_q}$ of degree bounded by ${\deg P< n(q-1)}$. Then 
${\sum_{\bfb \in \F^n}P(\bfb) =0}$.
\end{lemma}

\begin{proof}
This is Lemma~6.4 in~\cite{LN97}.
\end{proof}

\begin{lemma}
\label{lba}
Let~$\F$ be a finite field of odd characteristic and having more than~$3$ elements.  Further, let ${f(x,y), g(x,y)\in \F[x,y]}$ be quadratic forms over~$\F$. Then for ${c \in \F^\times}$ we have
$$
\sum_{\genfrac{}{}{0pt}{}{a,b\in \F}{g(a,b)=\pm c}}f(a,b)=0,
$$
where the sum is over the pairs ${(a,b)\in \F^2}$ such that ${g(a,b)=\pm c}$. 
\end{lemma}

\begin{proof}
Write ${q=|\F|}$, so that ${\F=\F_q}$. Then 
$$
\sum_{\genfrac{}{}{0pt}{}{a,b\in \F}{g(a,b)=\pm c}}f(a,b)= \sum_{a,b\in \F} f(a,b)(2-(g(a,b)-c)^{q-1}-(g(a,b)+c)^{q-1}). 
$$
We have 
$$
f(x,y)(2-(g(x,y)-c)^{q-1}-(g(x,y)+c)^{q-1})= -2f(x,y)g(x,y)^{q-1}+\big[\text{terms of degree $<2(q-1)$}\big],
$$
and Lemma~\ref{lln} implies that the sum from the assertion is equal to ${-2\sum_{a,b\in \F}f(x,y)g(x,y)^{q-1}}$.
The latter sum is 
${\sum_{\genfrac{}{}{0pt}{}{a,b\in \F}{g(a,b)\ne0}}f(a,b)}$,
which again by Lemma~\ref{lln} and by the assumption ${q>3}$ is equal to 
${-\sum_{\genfrac{}{}{0pt}{}{a,b\in \F}{g(a,b)=0}}f(a,b)}$.
If the quadratic form $g(x,y)$ is anisotropic over~$\F$ then the latter sum consists only of the term $f(0,0)$ and there is nothing to prove. And if it is isotropic then after a change of variables we may assume that ${g(x,y)=xy}$. Writing ${f(x,y)=\alpha x^2+\beta xy+\gamma y^2}$, the latter sum becomes
${(\alpha+\gamma)\sum_{a\in \F}a^2}$. 
Lemma~\ref{lln} implies that 
${\sum_{a\in \F}a^2=0}$
when~$\F$ has more than~$3$ elements. This completes the proof.
\end{proof}

\begin{proof}[Proof of Theorem~\ref{tu}]
Recall that  the orbit~$\OO$ consists of ${(x,y)\in \F_p^2}$ satisfying   ${\Xi x^2-y^2\in cH}$ with some ${c\in \F_p^\times}$.  Since ${H\ni -1}$, Lemma~\ref{lba} implies that the quadratic relations~\eqref{equadoo} hold true. Further, for each  ${c\in\F_p^\times}$ there is  exactly ${p+1}$ elements of $\F_{p^2}$ of norm~$c$, which implies that our orbit~$\OO$ has exactly ${(p+1)|H|}$ elements, and with our choice of~$m$ the divisibility conditions~\eqref{e12moo} hold true as well. Corollary~\ref{cecoun} now implies that ${u\in K(X_G)}$. 

Finally,~$u_\OO$ does not depend on the lifting. Indeed, if we choose two  different liftings respecting complex conjugation and obtain the products, say,~$u$ and~$u'$, then ${u/u'}$ is a $p$-th root of unity by Proposition~\ref{ptriv}. On the other hand, ${u,u'\in K(X_G)}$, which implies that ${u/u'\in K}$, a totally real field. Hence ${u=u'}$. The theorem is proved.
\end{proof}

\subsection{Galois action on the units}
\label{ssgalou}

Consider first the case of general algebraic curves. The proof of the following proposition is a standard exercise in Galois theory. 

\begin{proposition}
Let $K/k$ be a finite Galois extension of fields of characteristic~$0$, and let~$X$ be a projective curve defined (that is, having a geometrically irreducible model) over~$k$. Then the extension $K(X)/k(X)$ is Galois and the restriction map 
$$
\gal\bigl(K(X)/k(X)\bigr)\to \gal(K/k), \quad \sigma\mapsto\sigma\vert_K
$$ 
defines isomorphism of Galois groups. Further, for ${P\in X(k)}$ and ${u\in K(X)}$ we have ${u(P)\in K}$,  and given ${\sigma \in \gal\bigl(K(X)/k(X)\bigr)= \gal(K/k)}$ we have ${u^\sigma(P)=u(P)^\sigma}$. 
\end{proposition}

In our case the group
$$
\gal \bigl(K(X_G)/\Q(X_G)\bigr)=\gal(K/\Q)=G/G_H=\F_p^\times/H
$$
acts transitively and faithfully on the right $G_H$-orbits, and this action agrees with the Galois action: for ${\sigma \in \gal(K/\Q)=\F_p^\times/H}$ we have ${u_\OO^\sigma=u_{\OO\sigma}}$. Fixing an orbit~$\OO$ and putting ${U=u_\OO}$, we obtain the following.

\begin{proposition}
For ${P\in X_G(\Q)}$ we have ${U(P)\in K}$ and  ${U^\sigma(P)=U(P)^\sigma}$ for ${\sigma \in \gal(K/\Q)}$. 
\end{proposition}

Since distinct orbits are disjoint, Theorem~\ref{tmdr} and the discussion thereafter have the following consequence (recall that ${d=[K:\Q]=[\F_p^\times:H]}$). 

\begin{proposition}
\label{pmdrns}
The~$d$ principal divisors $(U^\sigma)$, ${\sigma \in \gal(K/\Q)}$, generate an abelian group of rank ${d-1}$, the only relation being ${\sum_\sigma(U^\sigma)=0}$. In particular, if ${d\ge 3}$ and ${\sigma \ne 1}$ then~$U$  and $U^\sigma$ are multiplicatively independent modulo the constants. 
\end{proposition}

Finally, equation~\eqref{eproduap} implies that
\begin{equation}
\label{eprodusig}
\prod_{\sigma \in \gal(K/\Q)}U^\sigma = \pm p^m
\end{equation}
Indeed, arguing as in the proof of Lemma~\ref{lprodofus}, we show that the left-hand side of~\eqref{eprodusig} is a rational constant. It also follows from our definitions that the left-hand side of~\eqref{eprodusig} is equal to 
$$
\prod_{\bfa\in M_P}g_\tilbfa^m.
$$
Razing this to $(12p/m)$th power, we obtain the left-hand side of~\eqref{eproduap}. This proves~\eqref{eprodusig}.

\section{The principal relation}
\label{sprinrel}

We retain the set-up of Section~\ref{scpu} and in particular that of Section~\ref{ssgalou}: 
\begin{itemize}

\item
${p\ge 5}$ is a prime number, $\zeta_p$ is a primitive $p$-th root of unity;

\item
${\bfa\mapsto \tilbfa}$ is a lifting of the set ${M_p=\F_p^2\smallsetminus\{(0,0)\}}$  which respects complex conjugation and satisfies~\eqref{e01};

\item
$G$ is the normalizer of a non-split Cartan subgroup of $\GL_2(\F_p)$, realized as in~\eqref{eexpl};

\item
$H$ is a subgroup of~$\F_p^\times$, ${H\ni -1}$;

\item
${m=2}$  or~${6}$ according to~\eqref{edefm}. 

\item
${K=\Q(\zeta_p)^H}$, ${d=[K:\Q]=[\F_p^\times:H]}$;

\item
$\OO$ is a fixed right $G_H$-orbit in $M_p$ and ${U=u_\OO}$ as defined in Theorem~\ref{tu}. It might be worth pointing out that 
\begin{equation}
\label{ecardo}
|H|= \frac{p-1}{d}, \qquad
|\OO|=(p+1)|H|=\frac{p^2-1}{d}.
\end{equation}

\end{itemize}

We fix a system ${\eta_1, \ldots, \eta_{d-1}}$  of fundamental units of the totally real field~$K$. We also  put
\begin{equation}
\label{eeta0}
\eta_0=\NN_{\Q(\zeta_p)/K}(1-\zeta_p).
\end{equation}
Clearly,~$\eta_0$ generates the prime ideal~$\gerp$ of~$K$ above~$p$; recall that ${\gerp^d=(p)}$.

Recall that we call a point ${P\in X_G(\Q)}$ \textsl{integral} if ${j(P)\in \Z}$. 
Proposition~\ref{psiu} implies  that for an integral point~$P$ on $X_G$, the principal ideal $\bigl(U(P)\bigr)$ is an integral ideal of the field~$K$, and moreover it is a power of~$\gerp$. Since ${\gerp^\phi=\gerp}$ for ${\phi\in \gal(K/\Q)}$, relation~\eqref{eprodusig} implies that ${\bigl(U(P)\bigr)= \gerp^m}$. 
Thus, we have 
\begin{equation}
\label{eup}
U(P) = \pm \eta_0^{b_0}\eta_1^{b_1}\cdots \eta_{d-1}^{b_{d-1}},
\end{equation} 
where ${b_0=m}$ and ${b_1, \cdots, b_{d-1}}$ are some rational integers depending on~$P$. 

The purpose of this section is to express the exponents~$b_k$ in terms of the point~$P$; more precisely, in terms of $q_c(P)$, where~$c$ is the nearest cusp to~$P$ (Section~\ref{ssnear}). This can be viewed as an analog of Equation~(20) on page~378 of~\cite{BH96}.

For ${\phi \in \gal(K/\Q)}$ we have\footnote{In the sequel we use the letter~$\phi$ rather than~$\sigma$ to denote elements of $\gal(K/\Q)$.}
$$
U^\phi(P) = \pm (\eta_0^\phi)^{b_0}(\eta_1^\phi)^{b_1}\cdots (\eta_{d_1}^\phi)^{b_{d-1}}.
$$
Fix an ordering on the elements of the Galois group: ${\gal(K/\Q)=\{\phi_0=\id, \phi_1, \ldots, \phi_{d-1}\}}$. 
Since the real algebraic numbers ${\eta_0,\eta_1, \ldots, \eta_{d-1}}$ are multiplicatively independent, the ${d\times d}$ real matrix 
$\bigl(\log |\eta_\ell^{\phi_k}|\bigr)_{0\le k,\ell\le d-1}$ is non-singular. Let $\bigl(\alpha_{k\ell}\bigr)_{0\le k,\ell\le d-1}$ be the inverse matrix. Then
\begin{equation}
\label{ebkprel}
b_k = \sum_{\ell=0}^{d-1}\alpha_{k\ell}\log|U^{\phi_\ell}(P)| \qquad (k=0, 1,\ldots, d-1). 
\end{equation}
We will call~\eqref{ebkprel} the \emph{principal relation}: It will play crucial role in our reduction and enumeration algorithms.

Combining~\eqref{ebkprel} with Proposition~\ref{passu}, we may express~$b_k$ in terms of $q_c(P)$. 
Let us introduce some notation.  Let~$c$ be a cusp of $X_G$.
Define the following quantities:
\begin{equation}
\label{edelvart}
\begin{aligned}
\delta_{c,k} &= -m\sum_{\ell=0}^{d-1}\alpha_{k\ell}\ell_{\OO\phi_\ell\sigma},&
%\varsigma_{c,\ell} = \prod_{\bfa \in \OO\phi_\ell\sigma_c}\varrho_\tilbfa^m, \qquad
\vartheta_{c,k}&= m\sum_{\ell=0}^{d-1}\alpha_{k\ell}\log|\varrho_{\OO\phi_\ell\sigma}| \qquad(k=0,1,\ldots, d-1);\\
\kappa&=\max_k\sum_{\ell=0}^{d-1} |\alpha_{k\ell}| ,&
 \Theta&=\kappa mp(p^2-1)d^{-1} ,
% \qquad\delta_{\max}=\delta_{\max,c}=\max_k|\delta_{c,k}|, \quad \vartheta_{\max}=\vartheta_{\max,c}=\max_k|\vartheta_{c,k}|. 
\end{aligned}
\end{equation}
where~$\ell_A$ and~$\varrho_A$ are defined in~\eqref{elara}, and~$\sigma$ is an element in $\Gamma(1)$ such that ${\sigma(i\infty)}$ represents~$c$. It follows from~\eqref{eordcoef} that $\delta_{c,k}$ and $\vartheta_{c,k}$ are independent of the choice of~$\sigma$.  

\begin{remark}
\label{rtriv}
It is easy  to see that ${\delta_{c,0}=0}$ and at least one of the numbers ${\delta_{c,1}, \delta_{c,2}, \ldots,\delta_{c,d-1}}$ is non-zero. Indeed, we have
\begin{equation}
\label{etriv}
\left(\begin{smallmatrix}\ord_cU^{\phi_0}\\\ord_cU^{\phi_1}\\ \vdots\\\ord_cU^{\phi_{d-1}}\end{smallmatrix}\right)=\bigl(\log |\eta_\ell^{\phi_k}|\bigr)_{0\le k,\ell\le d-1}\left(\begin{smallmatrix}\delta_{c,0}\\\delta_{c,1}\\ \vdots\\\delta_{c,d-1}\end{smallmatrix}\right).
\end{equation}
Multiplying both sides by the row vector $(1,\ldots,1)$ on the left, we obtain ${\delta_{c,0}=0}$. Further, since the column vector on the left of~\eqref{etriv} is non-zero, so is the column vector on the right. 
\end{remark}

\begin{proposition}
\label{pformula20}
Let~$P$ be an integral point on~$X_G$ and~$c$  its nearest cusp (that is, ${P\in \Omega_c}$). Then  for ${k=0, \ldots, d-1}$ we have
\begin{equation}
\label{erelation}
%\begin{aligned}
b_k = \delta_{c,k}\log|q_c(P)|^{-1}+ \vartheta_{c,k} \\
 %   &\hphantom{=} + m\sum_{\ell=0}^{d-1}\alpha_{k\ell}  \left( \sum_{\genfrac{}{}{0pt}{}{\bfa \in \OO\sigma_{\ell}\sigma_c}{0<\tila_1<1/2}}\log|1-q_c^{\tila_1}e^{2\pi i\tila_2}|+\sum_{\genfrac{}{}{0pt}{}{\bfa \in \OO\sigma_{\ell}\sigma_c}{1/2\le\tila_1<1}}\log|1-q_c^{1-\tila_1}e^{-2\pi i\tila_2}|\right) \\
%    &\hphantom{=} 
+O_1\bigl(\Theta |q_c(P)|^{1/p}\bigr). 
%\end{aligned}
\end{equation}
In particular, 
\begin{equation}
\label{ebq}
|b_k|\le |\delta_{c,k}| \log|q_c(P)|^{-1}+|\vartheta_{c,k}|+ \Theta. 
\end{equation}
\end{proposition}

\begin{proof}
Using~\eqref{esimpleestimateulog} and~\eqref{ecardo}, we obtain  
\begin{equation*}
\log|U(P)|= m\ell_{\OO\sigma} \log|q_c(P)|+m\log|\varrho_{\OO\sigma}|+O_1\bigl( 
mp(p^2-1)d^{-1}|q_c(P)|^{1/p}\bigr),  
\end{equation*}
and the same holds true when~$U$ is replaced by~$U^\phi$ and~$\OO$ by~$\OO\phi$ for any ${\phi\in \gal(K/\Q)}$. Now the result follows by~\eqref{ebkprel} and~\eqref{edelvart}. 
\end{proof}

%%%%%%%%%%%%%%%%%%%%%%%%%%%%%%%%%%%%%%%%%%%%%%%%%%%%%%%%%%%%%%%%%%%%%
%%%%%%%%%%%%%%%%%%%%%%%%%%%%%%%%%%%%%%%%%%%%%%%%%%%%%%%%%%%%%%%%%%%%%
\section{Outline of the algorithm} \label{secOutlineOfAlgorithm}
%%%%%%%%%%%%%%%%%%%%%%%%%%%%%%%%%%%%%%%%%%%%%%%%%%%%%%%%%%%%%%%%%%%%%
%%%%%%%%%%%%%%%%%%%%%%%%%%%%%%%%%%%%%%%%%%%%%%%%%%%%%%%%%%%%%%%%%%%%%

In order to compute the integral points on $X_\ns^+(p)$, we need to consider them on each $\calF(\sigma)$ (see Section~\ref{ssnear}), $\sigma\in\Sigma$, where $\Sigma$ is the full system of representatives of cosets $\Gamma\backslash\Gamma(1)$ as given in Section~\ref{sscusps}.
Let us fix one such~$\sigma$ until the end of Section~\ref{secEllipsoidsSieve}, and we denote by~$c$ the cusp represented by $\sigma(i\infty)$, so that ${\calF(\sigma)\subset\Omega_c}$.

The principal relation~\eqref{ebkprel} can be written as
\begin{equation}
\label{eq_bAlambda}
b = A\cdot \lambda,
\end{equation}
where $b=(b_k)_{0\leq k\leq d-1}$, $A = (\alpha_{k\ell})_{0\leq k,\ell\leq d-1}$, and $\lambda = (\log |U^{\phi_\ell}(P)|)_{0\leq \ell\leq d-1}$. 

For integral points $P$ with $j(P)\not\in\{1,2,\ldots,1727\}$  with closest cusp $c$, the associated $q_c$-parameter is real and non-zero, and more precisely 
$$
q_c(P)\in I_0 := [-e^{-\pi\sqrt3}, e^{-2\pi}] \smallsetminus\{ 0\}.
$$
Suppose further that $P$ lies in $\calF(\sigma)$, where we use the notation from Section~\ref{ssnear}.
The points in $\calF(\sigma)$ with real $q_c$-parameter are in bijective correspondence with $I_0$ via the map $q_c$.
Restricting $\lambda$ to those points of interest within $\calF(\sigma)$, we obtain a curve in $\RR^d$ parametrized over $I_0$, and we denote it $\lambda_\sigma(q_c): I_0\to\RR^d$ with ``variable'' $q_c\in I_0$.

Furthermore let us define a curve $\gamma_\sigma := A\cdot\lambda_\sigma$ in $\RR^d$, which is equally parametrized by~$q_c$ over the same domain~$I_0$, and which depends as well on the chosen $\sigma\in\Sigma$.
Recall that the matrix $A$ depends only on the choice of the fundamental units ${\eta_1,\ldots, \eta_{d-1}}$.
As ${b\in\ZZ^d}$, equation~\eqref{eq_bAlambda} tells us that each integral point $P$ in $\calF(\sigma)$ gives rise to an intersection of $\gamma_\sigma$ with the lattice~$\ZZ^d$.
Thus our algorithm will essentially (up to numerical issues) do the following.
\begin{enumerate}
\item For each $\sigma\in\Sigma$, compute all $q_c\in [-e^{-\pi\sqrt3}, e^{-2\pi}] \wo 0$ for which $\gamma_\sigma(q_c)\in\ZZ^d$, and compute the corresponding $j$-invariants.
\item For each such~$j$ and additionally all ${j\in\{1,2,\ldots,1727\}}$: If $j\in\ZZ$, compute the image type of an associated mod-$p$ Galois representation.
\item If furthermore the image of the representation is contained in the normalizer of a non-split Cartan subgroup, then output~$j$.
\end{enumerate}

\paragraph{An explicit bound for $|q_c|$  and first reductions.}
In Section~\ref{supperboundred} we  bound $|q_c|$ away from zero, and explain how this bound can be considerably improved using a classical reduction procedure going back to Baker and Davenport~\cite{BD69}.

\paragraph{Sieving for lattice points in the remaining parts of~$\gamma_\sigma$.}
After having applied the first reductions, we are left with considering two compact real intervals $I_+ \subset \RR_{>0}$ and $I_- \subset \RR_{<0}$ of $q_c$-parameters.
In Section~\ref{secEllipsoidsSieve} we will cover $\gamma_\sigma(I_+\cup I_-)$ with ellipsoids of small volume and then use the Fincke--Pohst algorithm in order to obtain only a few remaining candidates for $j(P)$, which then can be checked in the extra search.

\paragraph{Extra search.}
A small set of values $j(P)\in\ZZ$ can be checked separately by simply computing the image type of the corresponding mod-$p$ Galois representation, see Section~\ref{ssmissing}. We call this the extra search.

\medskip

Computing the initial bound for $|q_c|$ is very fast.
Starting from this bound, in principle we could use the sieve from Section~\ref{secEllipsoidsSieve} immediately without the first reductions from Section~\ref{supperboundred}. However the latter makes simpler estimates that require much less precision and are thus considerably faster.
When the first reduction method cannot continue further, the sieve using ellipsoids sieves away many candidates for~$j(P)$. Only at the end, for large $|q_c|$, the sieve eventually stops sieving away a lot of false candidates. At this point it becomes beneficial to simply use the extra search, as only a few values for $j\in\ZZ$ remain that correspond to those~$q_c$.

\begin{remark}[Taking care of numerical issues.]
\label{remRealIntervals}
Our algorithm needs to deal with real numbers. A curious and trivial
fact is that all real numbers in this paper (although $\RR$ is
uncountable) are exactly representable in a computer; however only
symbolically, not as floating point numbers or more generally as rational
numbers.
To do efficient computations we need to work with rational
approximations instead. Thus dealing with error estimates is
unavoidable.
To solve all arising numerical problems elegantly, we used throughout
\emph{interval arithmetic}, where a real number $x$ is replaced by an
interval $[a,b]$ with $a,b\in\QQ$, containing~$x$. This has the advantage
that at each step of the algorithm we have exact rational bounds for
the computed real numbers.
If at some point during runtime it turns out that these bounds are too
weak, we simply rerun the relevant parts of the program with larger
and larger precision, that is,  we represent the initial real numbers $x$
by shorter and shorter intervals $[a,b]$.
This will certainly increase the running time as the heights of $a,b$
will increase. But once the precision of the result is sufficient, we
are certain that the obtained bounds are correct.

Note also that whenever we found a lattice point on $\gamma_\sigma$, we
need to know which $q_c$-parameter it comes from in order to determine
the possible values of $j(P)$. Solving this equation to a high
precision can be very costly time-wise, which is why the first
reductions and the subsequent sieving steps are so important.
\end{remark}

\section{Baker's bound and its reduction}
\label{supperboundred}
In this section we bound  $q_c(P)$, using the bound for ${|j(P)|}$ obtained in~\cite{BS13} with Baker's method. Afterwards, we show how it can be improved using the  Baker--Davenport method~\cite{BD69}. 

\begin{theorem}
\label{thbasha}
Assume that ${p \ge 7}$. Let~$\delta$  be the smallest divisor of ${(p-1)/2}$ satisfying ${\delta \ge 3}$. Set 
$$
\mho= 30^{\delta+5} \delta^{-2\delta+4.5}p^{6\delta+5}(\log p)^2, \qquad \mho_0=\log\bigl(e^{\mho}+2079\bigr)
$$
Then
for any integral point~$P$  on $X_\ns^+(p)$
with nearest cusp~$c$ we have ${\log |j(P)| < \mho}$ and ${\log|q_c(P)^{-1}|\le \mho_0}$. 
\end{theorem}

\begin{proof}
The first  statement is a version of Theorem~1.1 of~\cite{BS13}. The second statement follows from the first one using~\eqref{e2079}.
\end{proof}

\begin{remark}
Modular units used in~\cite{BS13} are not ``economical'' (in the sens of this article). Using instead economical units would yield a slightly sharper bound than $\mho_0$. However, the quality of this bound is not of significant importance for the subsequent reduction process. Therefore we prefer to use the ``prêt à porter'' result from~\cite{BS13} rather than repeat the (rather technical) argument from that paper replacing the ``naive'' units used therein by economical ones.  
\end{remark}

We call~$\mho_0$  \emph{Baker's bound}. It is usually numerically huge (around $10^{100}$ for small $p$ and  even about $10^{1000}$ for $p=97$), and so are the implied bounds for the exponents $b_0, \ldots, b_{d-1}$, that can be obtained using~\eqref{ebq}; therefore they are not  suitable for direct enumeration of all possible vectors ${b = (b_0, \ldots, b_{d-1})}$. Moreover, checking whether such a candidate vector~$b$  comes from an integral point~$P$ is non-trivial and computationally expensive.

However, in practical situations the bound $\mho_0$ can be drastically reduced, using the numerical Diophantine approximations technique introduced by Baker and Davenport~\cite{BD69} and developed in~\cite{BH96,TW89} in the context of the Diophantine equation of Thue.

As in the previous section we fix a cusp~$c$ and consider integral points ${P\in \Omega_c}$. 
We shall usually omit the index~$c$, writing ${\delta_{c,k}=\delta_k}$ and ${\vartheta_{c,k}=\vartheta_{k}}$
for the quantities defined in~\eqref{edelvart}.  

%Now let us be more specific. 
As we have seen in Remark~\ref{rtriv}, at least one of the numbers ${\delta_1, \ldots, \delta_{d-1}}$ is non-zero.  To simplify notation we will assume 
that ${\delta_1\ne 0}$. 
We denote 
$$
B_0=|\delta_1| \mho_0+|\vartheta_1|+ \Theta, 
$$
so that ${|b_1|\le B_0}$ by~\eqref{ebq}.

Put
$$
\delta=\frac{\delta_{{2}}}{\delta_{{1}}}, \quad \lambda= \frac{\delta_{{2}}\vartheta_{{1}}-\delta_{{1}} \vartheta_{{2}}}{\delta_{{1}}}.
$$
Relation~\eqref{erelation} implies that
\begin{equation}
\label{eredula}
\left| b_{{2}}-\delta b_{{1}}+ \lambda \right| \le (1+|\delta|)\Theta |q_c(P)|^{1/p},
%\qquad\Theta=\kappa m(p+1)|H|. 
\end{equation}
which gives us up to~$\lambda$ a good rational approximation of~$\delta$.
The idea is to find by continued fraction expansion of~$\delta$ an integer~$r$ such that $r\delta$ is close to an integer but $r\lambda$ is not, which then gives us a better upper bound for ${\log|q_c(P)^{-1}|}$.
For this, we proceed as follows.
We fix a real number ${T\ge 2}$ (in our computations we initially take ${T=10}$).
Next, using continued fractions we find a ``good'' rational  approximation of $\delta$; precisely, we find a non-negative integer ${r \le T B_{0}}$ such that
$$
\| r \delta \| \le (T B_{0})^{-1}
$$
where $\|\cdot  \|$ is the distance to the nearest integer. 
Thus, $T$ controls the precision of this rational approximation on top of~$B_0$.
Now, if $r\lambda$ is not ``very close'' to the nearest integer (in practice if ${\| r \lambda \| \ge 2 T^{-1}}$) then we can bound ${|q_c(P)^{-1}|}$.
Indeed, 
multiply both sides of~\eqref{eredula} by~$r$.  
The right-hand side of the resulting inequality will be bounded from above by ${(1+|\delta|)\Theta TB_0 |q_c(P)|^{1/p}}$, and the  left-hand side %of the resulting inequality 
would be
$$
\left| rb_{{2}}-r\delta b_{{1}}+ r\lambda \right| \ge \|r\lambda\|-B_0\|r\delta\| \ge \|r\lambda\|-T^{-1}, 
$$
since ${|b_1|\le B_0}$.
This gives the following upper bound for  ${|q_c(P)^{-1}|}$,
\begin{equation}
\label{enbq}
\log|q_c(P)^{-1}|\le p\log\frac{(1+|\delta|)\Theta TB_0}{\|r\lambda\|-T^{-1}} =: \mho_1.
\end{equation}
In the case when ${\|r\lambda\|<2T^{-1}}$ we increase~$T$ (in our computations we replace it by $10T$) and restart, until ${\|r\lambda\|\geq 2T^{-1}}$.
In practise only very few such iterations are needed: Heuristically, $\|r\lambda\|$ can be thought of a random number in $[0,1]$ (with respect to the Lebesgue measure), and thus a particular $T$ works with `probability'~$1-2T^{-1}$.
As $\log(TB_0)$ is almost $\log(B_0)$, the precise value of $T$ is not that important.

Since~$\mho_1$ depends logarithmically on~$\mho_0$, it is expected to be much smaller  than~$\mho_0$, and in practice it is.

We then repeat the same procedure, but this time with~$\mho_1$ instead of~$\mho_0$, obtaining for $\log|q_c(P)^{-1}|$ a new reduced bound~$\mho_2$, and so on. 
We stop this reduction process once it does not improve the previous bound on $\log|q_c(P)^{-1}|$ by more than~$1\%$.
In practice, three to four iterations of this procedure suffice for that.
We call~$\hatmho$ the obtained reduced bound for
$\log|q_c(P)^{-1}|$. %  and~$\hatB$, respectively. 
In practice~$\hatmho$ is around $200$ for small $p$ and about $2200$ for $p=97$.

\begin{remark}[Geometry intuition behind the reduction process]
Recall from Section~\ref{secOutlineOfAlgorithm} that we are interested in the real $q_c$-para\-me\-ters $q_c\in I_0$ at which $\gamma_\sigma(q_c)$ is a lattice point $b\in\ZZ^d$.
For $|q_c|\to 0$ the $\ell_2$-norm of $\gamma_\sigma$ tends to infinity, and the asymptotic direction of $\gamma_\sigma$ is given by the vector $(\delta_{0},\ldots,\delta_{d-1})$, compare with Remark~\ref{rtriv}.
Moreover, the Baker bound $\mho_0$ for $\log|q_c^{-1}|$ together with~\eqref{ebq}  restricts our search to a certain hypercube   in $\RR^d$.
In informal, more geometric terms, the reduction procedure can be described as follows. 
We project~$\gamma_\sigma$ to its first two coordinates, and we try to find a large part of the domain of~$\gamma$ where this projection does not intersect the lattice points $\ZZ^2\subset\RR^2$.
For this we approximate the slope of the projected $\gamma_\sigma$ by a rational number; in other words we make a change of coordinates via a matrix in $\GL_2(\ZZ)$, sending $\ZZ^2$ to $\ZZ^2$, such that the projected $\gamma_\sigma$ becomes asymptotically almost horizontal. That is, the vertical coordinate changes asymptotically very slowly, which means that it stays over a large part of the domain between two integers.
Then using the error bounds of~\eqref{erelation} we deduce a new upper bound $\mho_1$ for $\log|q_c^{-1}|$.
\end{remark}

\section{Sieving lattice points on  \texorpdfstring{$\gamma_\sigma$}{gs}}
\label{secEllipsoidsSieve}

\subsection{Ellipsoid sieve}

We continue with the set-up of Section~\ref{secOutlineOfAlgorithm}.
For notational simplicity, we don't keep $\sigma\in\Sigma$ and $c$ in the notation. That is, we write $\lambda$ for $\lambda_\sigma$, $\gamma$ for $\gamma_\sigma$, and $q$ for $q_c$.
Suppose that we want to find all $q$ in a certain interval within the domain of $\gamma$ at which $\gamma(q)\in\ZZ^d$.

One serious numerical problem is, that we can compute the modular units $U$, and thus $\gamma$, only up to an arbitrary finite precision, but never exactly.
Moreover, a better precision (i.e. a better error bound) requires longer running time and more memory consumption.
Thus we are looking for ways to algorithmically bound $\gamma$ in such a way that only very few candidates~$j$ remain.
For this we use two different methods:
\begin{enumerate}
\item We will cover (the image of) $\gamma$ with exactly computed ellipsoids, such that each of the ellipsoids intersects~$\ZZ^d$ in very few points, which can be determined using an enhanced version of the Fincke--Pohst algorithm~\cite{FP85}. 
\item If a subinterval $I$ of the domain of $\gamma$ has the property, that the restriction $\gamma|_I$ has a coordinate with provably positive or negative derivative, we can use the bisection method or Brent's method~\cite{Br71} to compute all $q$ for which this coordinate is integral. 
\end{enumerate}

By default we try to use the first of the two methods, because it is considerably faster in most situations.
Only at places where $\gamma$ intersects $\ZZ^d$ or is very close to doing so, the second method becomes preferable.
Next we discuss both these two methods in detail.

\paragraph{Covering $\gamma(I)$ with one ellipsoid.}
Let ${I=[q_1,q_2] \subset\RR\smallsetminus\{0\}}$ be an interval in the domain of~$\gamma$. Recall that
${U=u_\OO}$ is a product of Siegel functions $g_a$ (Theorem~\ref{tu}), whose leading terms come from the factors in front of the infinite products in~\eqref{epga}.
Therefore the leading term of $\log|U|$ is affine linear in~$\tau$. 

The remainder of $\log |U|$ is an infinite sum, which can be estimated by a finite sum plus an error term, see Corollary~\ref{cappru} in the Appendix. 
We use real and complex interval arithmetic, using that~$q$ lies in the given interval~$I$, and obtain constant lower and upper bounds for the remainder of $\log |U|$ that hold for all $\tau$ in the interval given by $q\in I$.

As the coordinates of~$\lambda$ are sums of terms $\log |U|$, it follows that the image of $\lambda|_I$ lies in the Minkowski sum of a line segment (coming from the leading terms of the terms $\log |U|$) and an axis-parallel cube (coming from the remaining terms of the summands $\log |U|$).
Let us call the line segment $S$, and the cube~$C$.
Thus, $\lambda(I)\subseteq S+C$, where ``$+$'' denotes the Minkowski sum.
We look for an ellipsoid $E$ of small volume that contains $S+C\supseteq \lambda(I)$.
For this we need some preparation.

Let $Q$ be a symmetric positive definite $d\times d$ matrix, which gives rise to a quadratic form on $\RR^d$.
In the following we write $E_Q := \{x\in\RR^d\st x^t Q x \leq 1\}$ for its unit ball, which is a euclidean ellipsoid.

\begin{lemma}[Product of ellipsoids]
\label{lemProductOfEllipsoids}
Let $Q_1$, $Q_2$ be positive definite matrices giving rise to quadratic forms on $\RR^{d_1}$ and $\RR^{d_2}$, respectively.
Then the ellipsoid $E_Q$ of smallest volume containing $Q_1\times Q_2$ is given by the block-diagonal matrix $Q$ with the two blocks $\frac{d_1}{d_1+d_2}Q_1$ and $\frac{d_2}{d_1+d_2}Q_2$.
\end{lemma}

\begin{lemma}[Projecting ellipsoids]
\label{lemProjectingEllipsoid}
Let $Q$ be a positive definite $d\times d$ matrix.
Let $P$ be a $d'\times d$ matrix of rank $d'$ representing a linear surjection $\RR^d\to\RR^{d'}$.
The image of $E_Q$ under $P$ is given by $P(E_Q) = E_{Q'}$ with $Q' = (P Q^{-1} P^t)^{-1}$.
\end{lemma}

\begin{proposition}[Convex hulls of ellipsoids]
\label{propConvexHullOfEllipsoids}
Let $Q$ be a positive definite $d\times d$ matrix, and let $a\in\RR^d$ be a vector.
Define $E_{Q'}$ as the ellipsoid given by
$Q' = \frac{1}{d+1} (\frac{1}{d}Q^{-1} + a^t a)^{-1}$.
Then $E_{Q'}$ is an ellipsoid that contains both translates $E_Q+a$ and $E_Q-a$ of $E_Q$.
\end{proposition}

\begin{proof}
Let $P$ be the $d\times (d+1)$ matrix $(\id_d | a)$, i.e.~a $d\times d$ identity matrix with an augmented column given by vector~$a$.
The convex hull of $E_Q+a$ and $E_Q-a$ can be written as $P(E_Q\times [-1,+1])$. Note that $[-1,+1] = E_{\id_1}$, where $\id_1$ is the $1\times 1$ identity matrix.
Using Lemma~\ref{lemProductOfEllipsoids} we find an ellipsoid $E_{\wt Q}$ containing $E_Q\times [-1,+1]$, namely the one given by a block matrix $\wt Q$ with a $d\times d$ block $\frac{d}{d+1}Q$ and a $1\times 1$ block with entry $\frac{1}{d+1}$.
Define $Q'$ as in Lemma~\ref{lemProjectingEllipsoid} such that $P(E_{\wt Q}) = E_{Q'}$.
Then one quickly checks that $Q'$ is given as in the assertion of the proposition, i.e. $Q' = \frac{1}{d+1} (\frac{1}{d}Q^{-1} + a^t a)^{-1}$, and we have
\[
E_Q \pm a \subseteq \conv(E_Q-a,E_Q+a) = P(E_Q\times [-1,+1]) \subseteq P(E_{\wt Q}) = E_{Q'}.
\]
\end{proof}

Now let us return to the problem of finding an ellipsoid $E$ containing $S+C\supseteq \lambda(I)$.
It is easy to write down the ellipsoid $E_C$ of smallest volume that contains~$C$; one can obtain it by iterating Lemma~\ref{lemProductOfEllipsoids}, using that $C$ is a product of $1$-dimensional ellipsoids.
By translating the coordinate system, we may assume that the segment $S$ is centered at the origin, its endpoints being $\pm a$.
Then Proposition~\ref{propConvexHullOfEllipsoids} yields an ellipsoid that contains both $E_C\pm a$, and thus $S+C$, and thus also $\lambda(I)$.

%\begin{lemma}[Affine transformation of ellipsoids]
%Let $Q$ be a positive definite $d\times d$ matrix giving rise to an ellipsoid $E_Q$.
%Let $A$ an invertible $d\times d$ matrix giving rise to a linear automorphism of~$\RR^d$.
%Then $A(E_Q) = E_{Q'}$ where $Q' = A^{-t} Q A^{-1}$.
%\end{lemma}
%
%\begin{proof}
%\[
%A(E_Q) = \{Ax\in\RR^d\st x^t Q x \leq 1\}
%= \{x\in\RR^d\st x^t A^{-t} Q A^{-1} x \leq 1\}
%= E_{Q'}.
%\]
%\end{proof}

As $\gamma(I) = A\cdot \lambda(I)$ is just an affine image of $\lambda(I)$, with Lemma~\ref{lemProjectingEllipsoid} we obtain immediately an ellipsoid containing $\gamma(I)$. Let us call this ellipsoid $E_I$.

\paragraph{Computing lattice points on $\gamma(I)$.}
We are interested in the $q$-parameters such that $\gamma(q)\in\ZZ^d$.

For this, in the previous section we computed an ellipsoid $E_I\supseteq\gamma(I)$.
With the Fincke--Pohst algorithm, one can compute all lattice points in an ellipsoid $E_I$.
We note that one needs to adjust the original Fincke--Pohst algorithm in two ways.
\begin{enumerate}
\item
As the dimensions $d$ become relatively high (say larger than 8), it is in practice necessary to first LLL-reduce the basis of the lattice with respect to the quadratic form that defines the ellipsoid.
Otherwise the Fincke--Pohst algorithm can become too slow.

\item
Our ellipsoid $E_I$ is in general not centered at a lattice point; whereas the original Fincke--Pohst algorithm is for ellipsoids centered at the origin. This generalization is indeed not difficult to implement (see for example our source code for our precise implementation).
\end{enumerate}

Note that everywhere in the code we use interval arithmetic.
That is, our computations do not compute the exact real numbers, but instead exact intervals in which the correct results of the computations lie; compare with Remark~\ref{remRealIntervals}.
In the Fincke--Pohst algorithm this ensures that we get indeed all lattice points in $E_I$, and possibly a few more ``false candidates''.
The higher the underlying precision is that we are using, the fewer false candidates there will be. 

\paragraph{Examining lattice points on $\gamma(I)$.}
Having computed the lattice points $E_I\cap\ZZ^d$ (plus possibly some false candidates), there are basically two possibilities:

\begin{enumerate}
\item
It may happen that Fincke--Pohst returns no lattice point.
Then we have a proof that $E_I\cap\ZZ^d$ is empty, and hence $\gamma(I)\cap\ZZ^d$ is empty.
Thus $X^+_\ns(p)$ contains no integral point within the hyperbolic triangle $\calF(\sigma)$ with $q$-parameter in~$I$.

\item It may happen that Fincke--Pohst returns at least one lattice point.
\label{ireturns}

\begin{enumerate}
\item
\label{iitwopoints}
One way to continue is to split $I$ into two smaller intervals $I=I_1\cap I_2$ and continue with $I_1$ and $I_2$ recursively.
The intuition is that the corresponding ellipsoids $E_{I_1}$ and $E_{I_2}$ should be of much smaller volume and should thus contain fewer lattice points. 
In the current implementation we do exactly this when Fincke--Pohst returns at least two lattice points, and we split $I$ into two equal pieces $I_1$ and $I_2$ in the logarithmic scale (which is a good splitting point in practise; it corresponds to bisecting the interval in $\tau$-coordinates).

\item
When Fincke--Pohst returns exactly one lattice point $v$, it may happen that it lies actually on $\gamma(I)$ and which we will not get rid of by splitting $I$ into smaller pieces.
Thus we try to compute which $j$-invariant(s) it corresponds to.
This is non-trivial for numerical reasons.
Our implementation does the following.
We compute the derivative of $\gamma|_I$ in interval arithmetic.
If in this way we cannot prove that at least one coordinate of $\gamma|_I$ has everywhere positive or everywhere negative derivative, we continue as in~\ref{iitwopoints} by splitting~$I$ into two pieces and go deeper into the recursion.
Otherwise we know that one coordinate, say the $k$'th one, is strictly monotone.
In particular, $\gamma|_I$ may go through the obtained lattice point $v$ at most once.
Then we use the bisection method (as always in interval arithmetic) to find the $q$-parameter (well, only an interval $I_q$ that contains this $q$-parameter) such that the $k$'th coordinates satisfy $\gamma_k(q) = v_k$.
It may happen that $\gamma_k$ is larger (or smaller) than $v_k$ at both endpoints of $I$, in which case we also proved that $\gamma(I)\cap\ZZ^d$ is empty.
Having found $I_q$, we compute in interval artihmetic $\gamma(I_q)$, which is a vector whose coordinates are all real intervals, i.e. $\gamma(I_q)$ is a cube. We check whether it contains the obtained lattice point~$v$.
If not, again we know that $\gamma(I)\cap\ZZ^d$ is empty.
Otherwise, we compute the set of all integral values for $j$ that correspond to $q$-parameters in $I_q$.
If there are no such integral values for $j$, we are of course again finished.
If there is only one such value for $j$, we test it in the extra search, see Section~\ref{ssmissing}.
However if there are more than one such value, we increase the precision of our computations, and rerun the bisection method that computes the real interval $I_q$, until at some point $I_q$ is small enough.
\end{enumerate}

\end{enumerate}

\begin{remark}
In an earlier version of this paper (see the first preprint version on the arXiv) we did not cover $\gamma$ by ellipsoids.
Instead, in the language of this section, we concentrated on one of the coordinates $\gamma_i$ of $\gamma$ and solved $\gamma_i(q) = k$ for all integral values of~$k$ within the range that is left after the first reductions of Section~\ref{supperboundred}. 
The ellipsoid method has two advantages:
Firstly, for~$q$ close to zero it is faster, as a single ellipsoid can cover many choices for~$k$.
In principle one could even replace the first reductions from Section~\ref{supperboundred} by suitable ellipsoids\footnote{This is analogous to the situation for $S$-unit and Mordell equations, where one can replace the classical first reduction of the initial height bound by ellipsoids~\cite{vKM16}.}, however the above reduction is much faster as it requires little precision.
Secondly, solving $\gamma_i(q) = k$ for~$q$ runs into considerable numerical issues, especially for larger $|q|$ due to slow convergence.
This being said, if done right, a coordinate-wise method might be a way to improve on our algorithm in the future, in particular as a middle step between the ellipsoid approach and the extra search, where the largest running time improvements are possible.
\end{remark}

\subsection{Extra search}
\label{ssmissing}

Assume that it only remains to verify for a few values for $j(P)\in\ZZ$, whether they come from an integral point~$P$.
In practice, these values will contain the numbers $1,2,\ldots,1727$, which come from potential integral points with non-real $q_c$-parameter. Moreover it will contain the values of $j$ that could not be excluded during the sieve of Section~\ref{secEllipsoidsSieve}, as well as all ${j\in\ZZ}$ with ${|j|\leq j_0}$ for some $j_0$ such as ${j_0=2^{16}}$ (or even larger for big~$p$), as for small~$j$ the sieve is slower than the  direct extra search described below. 
Also, the sieve has to be performed for each ${\sigma\in\Sigma}$, and the set~$\Sigma$ is of cardinality ${p(p-1)/2}$, quite a big number for big~$p$. 

Recall that to an elliptic curve~$E/\Q$ and a prime number~$p$ we associate a Galois representation
$$
\rho_{E,p}:\gal_\Q\to \GL(E[p])\cong \GL_2(\F_p), 
$$
which is defined by the natural action of the absolute Galois group~$G_\Q$ on the torsion group $E[p]$. Points in $X_\ns(p)(\Q)$ correspond  to the elliptic curves $E/\Q$ such that the image of $\rho_{E,p}$ is contained in the normalizer of a non-split Cartan subgroup of $\GL_2(\F_p)$. 

It is known that, if this latter property holds for some elliptic curve $E/\Q$ with ${j(E)\ne 0,1728}$, then it holds for any quadratic twist of~$E$, that is, for any other elliptic curve~$E'$ with ${j(E')=j(E)}$. Indeed,~$E'$ is isomorphic to~$E$ over some field~$K$ of degree at most~$2$. Denote by~$\chi_K$ the character of~$G_\Q$ corresponding to~$K$. Then  ${\rho_{E',p}=\rho_{E,p}\chi_K}$.  Hence if the image of $\rho_{E,p}$  is contained in the normalizer of a non-split Cartan subgroup, then so is the image of~$\rho_{E',p}$. 

Hence, if we fix ${j\in \Q}$, distinct from~$0$ and~$1728$, then,  to verify whether $X_\ns^+(p)$ has a rational point~$P$ with  ${j(P)=j}$, it suffices to verify for at least one curve $E/\Q$ with ${j(E)=j}$ whether the image of $\rho_{E,p}$ is contained in the normalizer of a non-split Cartan subgroup. This can be readily accomplished within \SAGE~\cite{sage} using the functions \emph{E $=$ EllipticCurve\_from\_j(j)} and
 \emph{E.galois\_representation().image\_type(p)}.

\subsection{Running time of the algorithm}
\label{secRunningTime}
%The following tables lists the running times that our algorithm took for the primes $7\leq p \leq 97$.

In order to obtain Theorem~\ref{thmMain}, we let our algorithm run for all primes $7\leq p\leq 97$ in parallel on the PlaFRIM computer cluster.
%The following table lists $t$, which is the running time in seconds multiplied by the number of CPUs that were in use.
For any such $p$, let $t$ denote the running time in seconds multiplied by the number of CPUs that were in use; see Table~\ref{tabRunningTimes}.
Thus, $t$ is an upper bound for the time in seconds that our algorithm would take for a prime $p$ on a single CPU.
%Let $t$ denote the running time in seconds that our algorithm would take for a particular prime $p$ on a single CPU. 
The computation for $p=79$ had the highest running time of almost $3.5$~CPU years.
The total running time was almost $16$ %$15.8$ 
CPU years.

\newcommand{\ruleToggle}[1]{#1}
\begin{table}[htb]
%\small
\begin{center}
\footnotesize
\begin{tabular}{ccc}
\ruleToggle{\toprule}
$p$ & $t$ & $t/p^4$ \\
\midrule%\hline
7 & $6.7 \cdot 10^{3}$ & $2.76$ \\
11 & $6.8 \cdot 10^{3}$ & $0.47$ \\
13 & $2.9 \cdot 10^{4}$ & $1.01$ \\
17 & $2.8 \cdot 10^{4}$ & $0.33$ \\
19 & $1.2 \cdot 10^{5}$ & $0.89$ \\
23 & $9.0 \cdot 10^{4}$ & $0.33$ \\
\ruleToggle{\bottomrule}
\end{tabular}
\quad
\begin{tabular}{ccc}
\ruleToggle{\toprule}
$p$ & $t$ & $t/p^4$ \\
\midrule%\hline
29 & $2.4 \cdot 10^{5}$ & $0.34$ \\
31 & $8.3 \cdot 10^{5}$ & $0.90$ \\
37 & $1.8 \cdot 10^{6}$ & $0.94$ \\
41 & $1.2 \cdot 10^{6}$ & $0.41$ \\
43 & $2.7 \cdot 10^{6}$ & $0.77$ \\
47 & $2.0 \cdot 10^{6}$ & $0.40$ \\
\ruleToggle{\bottomrule}
\end{tabular}
\quad
\begin{tabular}{ccc}
\ruleToggle{\toprule}
$p$ & $t$ & $t/p^4$ \\
\midrule%\hline
53 & $4.3 \cdot 10^{6}$ & $0.55$ \\
59 & $1.3 \cdot 10^{7}$ & $1.03$ \\
61 & $1.8 \cdot 10^{7}$ & $1.24$ \\
67 & $3.0 \cdot 10^{7}$ & $1.47$ \\
71 & $2.6 \cdot 10^{7}$ & $1.01$ \\
73 & $6.9 \cdot 10^{7}$ & $2.42$ \\
\ruleToggle{\bottomrule}
\end{tabular}
\quad
\begin{tabular}{ccc}
\ruleToggle{\toprule}
$p$ & $t$ & $t/p^4$ \\
\midrule%\hline
79 & $1.2 \cdot 10^{8}$ & $2.83$ \\
83 & $6.1 \cdot 10^{7}$ & $1.28$ \\
89 & $1.0 \cdot 10^{8}$ & $1.59$ \\
97 & $6.2 \cdot 10^{7}$ & $0.70$ \\
\ruleToggle{\bottomrule}
~\\
~\\
\end{tabular}
\end{center}
\caption{Running times $t$ in CPU seconds for computing $X_\ns^+(p)(\ZZ)$.}
\label{tabRunningTimes}
\end{table}

For each $p\neq 97$ the computation was based on the field $K={\Q(\zeta_p+\bar\zeta_p)}$. 
For $p=97$ we took for $K$ instead the subfield of degree $16$, as during a test-run on a single hyperbolic triangle $\calF(\sigma)$ this choice improved the running time by a factor of~$4.5$.
The running time for $p<100$ in our computations seems to grow slightly faster than~$p^4$, but as for the asymptotic complexity for our algorithm this is a strict underestimate.
A $p^2$ factor already comes from the number of hyperbolic triangles $\calF(\sigma)$ that cover $X^+_\ns(\CC)$ over which the algorithm iterates.
For each $\calF(\sigma)$, the bottleneck (at least in the current range) is the computation of integral points in the ellipsoids.
For each such ellipsoid, the Fincke--Pohst algorithm needs to be applied to a LLL-reduced quadratic form (otherwise it becomes impractical for moderately large $d$), and it also needs $O(d^3)$ time to compute the Fincke--Pohst form of the ellipsoid.
Moreover the number of covering ellipsoids and their required numerical precision is hard to estimate a priori, which makes even an upper bound for the running time of our algorithm currently out of reach.

Note that we did not prove that the algorithm actually terminates.
Something very special would need to happen for it not to terminate, such as the existence of an integral point out of the reach of the extra search at which the derivative of $\gamma$ vanishes.
In any case, if the algorithm will not terminate for some larger $p$ in the future, presumably there will be a quick fix of the program.

If the algorithm terminates then the output is proved to be correct, i.e. the obtained set of integral points on $X^+_\ns(p)$ is provably complete.

%TODO: 

%Y4.2: It must be obvious but I am really unclear on the need at the end for splitting the computations to find the nearest cusp, because at the end + of the method for X ns (p) (because the considered modular units are Galois conjugates, and the cusps as well), the situation seems perfectly symmetric to me, so I suggest adding a remark on why it is not the case.
%---> Answer: Numerics (better convergence), but where to answer that? 

%Y4.3: I guess using ellipsoids is the right way to go, but why is it not efficient to do an analysis coordinate by coordinate ? Is it because one needs to compute A with great precision? Could it work as a first sieve before Fincke-Pohst? In any case, I think it needs to be addressed.

\appendix

\section{Approximate formulas for Siegel functions, modular units and their derivatives}
\label{appendix}

In this appendix we collect formulas we were using  in the computation of Section~\ref{secEllipsoidsSieve} to numerically approximate modular units and their derivatives. 

We start by approximating Siegel functions. We use notation of Subsection~\ref{ssklesie}. 

\begin{proposition}
\label{papprsiegel}
Let~$\nmax$ be a positive integer. Then for a non-zero ${\bfa=(a_1,a_2) \in \Q^2\cap[0,1)^2}$ and ${\tau\in \HH}$ we have
\begin{equation}
\label{eqLogga}
\begin{aligned}
\log |g_\bfa(\tau)|&= \ell_\bfa \log|q|+ \log|\varrho_\bfa|
\\& \hphantom{=} +  
\sum_{\genfrac{}{}{0pt}{}{n=0}{(n,a_1)\neq (0,0)}}^{\nmax-1} \log |1-q^{n+a_1}e^{2\pi i a_2}|
+
\sum_{n=0}^{\nmax-1} \log |1-q^{n+1-a_1}e^{-2\pi i a_2}|
\\& \hphantom{=} +
O_1\left(  \frac{|q|^{\nmax+a_1}+|q|^{\nmax+1-a_1}}{(1-|q|)^2}\right).
\end{aligned}
\end{equation}
\end{proposition}

\begin{remark}
The term $\log|\varrho_\bfa|$ vanishes unless ${a_1=0}$. In this latter case it might have been included into the first sum by omitting the condition ${(n,a_1)=(0,0)}$. However, we prefer to write~\eqref{eqLogga} as we do, because we want to separate the constant and the non-constant terms. 
\end{remark}

\begin{proof}
The proof is similar to that of Proposition~\ref{pellaass}. Using the inequality  ${\bigl|\log (1+z)\bigr|  \le |z|/(1-|z|)}$ for ${|z|<1}$, we estimate, for ${n\ge 1}$, the terms of the product expansion~\eqref{epga}  as
\[
\bigl|\log(1-q^{n+\tila_1}e^{2\pi i\tila_2})\bigr|\le   \frac{|q|^{n+\tila_1}}{1-|q|^n}, \qquad 
\bigl|\log(1-q^{n+1-\tila_1}e^{-2\pi i\tila_2})\bigr|\le\frac{|q|^{n+1-\tila_1}}{1-|q|^n}. 
\]
Adding this up for all ${n\geq \nmax}$, we bound the error term in~\eqref{eqLogga} as 
$$
\sum_{n=\nmax}^\infty \frac{|q|^{n+\tila_1}+|q|^{n+1-\tila_1}}{1-|q|} 
%= \frac{|q_c|^{\nmax+\tila_1}+|q_c|^{\nmax+1-\tila_1}}{1-|q_c|} \cdot \sum_{k=0}^\infty |q_c|^k
= \frac{|q|^{\nmax+\tila_1}+|q|^{\nmax+1-\tila_1}}{(1-|q|)^2}. 
$$
This proves the proposition. 
\end{proof}

The absolute value $|g_\tilbfa(\tau)|$ is $1$-periodic as a function of~$\tau$. Hence it can be viewed as a function of ${q=e^{2\pi i\tau}}$. We  need an estimate for the logarithmic $q$-derivative $\frac{d}{dq}\log |g_\tilbfa|$. 

\begin{proposition}
\label{pderivsieg}
Let~$\tilbfa$,~$\tau$ and~$\nmax$ be as in Proposition~\ref{papprsiegel}. Assume that ${\tau\in \ocalF}$. Then
\begin{equation}
\label{eqDerivedLogga}
\begin{aligned}
\left(\frac{d}{dq}\log |g_\tilbfa|\right)(\tau)&= \ell_\tilbfa \re\frac{1}{q}+\sum_{n=0}^{\nmax-1}
\re \left(
%\log |1-q^{n+a_1}e^{2\pi i a_2}|
\frac{-(n+\tila_1)q^{n-1+\tila_1}e^{2\pi i\tila_2}}{1-q^{n+\tila_1}e^{2\pi i \tila_2}}
+
%\log |1-q^{n+1-a_1}e^{-2\pi i a_2}|
\frac{-(n+1-\tila_1)q^{n-\tila_1}e^{-2\pi i\tila_2}}{1-q^{n+1-\tila_1}e^{-2\pi i \tila_2}}
\right)
\\& \hphantom{=} +
O_1\left(
%\bigg(
%\frac{1}{(1-|q|)^2}-\sum_{n=1}^{\nmax-1}n|q|^{n-1}
%\bigg)
\frac{\nmax|q|^{\nmax-1}-(\nmax-1)|q|^{\nmax}}{(1-|q|)^2}
\left(\frac{1}{1-|q|^{\nmax+\tila_1}}+\frac{1}{1-|q|^{\nmax+1-\tila_1}}\right)
\right).
\end{aligned}
\end{equation}
\end{proposition}

\paragraph{Proof}
For any smooth function ${f:\R\to\C^\times}$, we have $\frac{d}{dx}\log|f(x)| = \re\frac{f'(x)}{f(x)}$.
Using the  product expansion~\eqref{epga}, this already explains the main term in~\eqref{eqDerivedLogga}. 
For the $O_1$-term we bound
$$
\re\left(\frac{-(n+\tila_1)q^{n-1+\tila_1}e^{2\pi i\tila_2}}{1-q^{n+\tila_1}e^{2\pi i \tila_2}}\right)
\leq
\left|\frac{-(n+\tila_1)q^{n-1+\tila_1}e^{2\pi i\tila_2}}{1-q^{n+\tila_1}e^{2\pi i \tila_2}}\right|
\leq
\frac{(n+\tila_1)|q|^{n-1+\tila_1}}{1-|q|^{n+\tila_1}}
$$
For the denominator, we simply bound ${1-|q|^{n+\tila_1} \geq 1-|q|^{\nmax+\tila_1}}$ for ${n\geq \nmax}$.
We want to bound the numerator for $n\geq \nmax$ as
$$
(n+\tila_1)|q|^{n-1+\tila_1}\leq n|q|^{n-1}.
$$
This follows from the fact that the function ${x\mapsto x |q|^{x-1}}$ is decreasing for ${x\ge 1}$ whenever ${|q| < e^{-1/x}}$, which is true because ${\tau\in \calF}$, %(because its derivative is $(1+\wt n\log|q_c|)|q_c|^{\wt n-1}$), 
and then plugging into this function the two arguments ${x = n}$ and ${x = n+\tila_1}$.

Thus, for ${n\geq \nmax}$, we obtain
$$
\sum_{n=\nmax}^\infty (n+\tila_1)|q|^{n-1+\tila_1}
\leq 
\sum_{n=\nmax}^\infty n|q|^{n-1}
%= 
%\sum_{n=0}^\infty n|q_c|^{n-1} - \sum_{n=0}^{\nmax-1} n|q_c|^{n-1}
=
%\frac{1}{(1-|q_c|)^2} - \sum_{n=0}^{\nmax-1} n|q_c|^{n-1},
\frac{\nmax|q|^{\nmax-1}-(\nmax-1)|q|^{\nmax}}{(1-|q|)^2}. 
$$
%as $\sum_{n=0}^\infty n|q_c|^{n-1} = \frac{1}{(1-|q_c|)^2}$ and $\sum_{n=0}^{\nmax-1} n|q_c|^{n-1} = \frac{(\nmax-1)|q_c|^\nmax - \nmax |q_c|^{\nmax-1}+ 1}{(1-|q_c|)^2}$, where the last equation is obtained by deriving the geometric series $\sum_{n=\nmax}^\infty |q_c|^n = \frac{|q_c|^\nmax}{1-|q_c|}$.
This explains the first half of the $O_1$-term. The other half for exponent ${n+1-a_1}$ instead of ${n+a_1}$ is treated in the same way.
\qed

\bigskip

Now it is easy to obtain a similar result for modular units. 

\begin{corollary}
\label{cappru}
Let~$\nmax$ be a positive integer. Then in the set-up of Proposition~\ref{passu}, for every ${P\in \ocalF(\sigma)}$ we have
\begin{equation*}
%\label{eqLogUOofP}
\begin{aligned}
\log |u(P)|&= \frac{\ord_cu}{p} \log|q_c(P)|+ \log|\varsigma_c|
\\& \hphantom{=} +  
m\sum_{n=0}^{\nmax-1}
\left(
\sum_{\genfrac{}{}{0pt}{}{\bfa \in \OO\sigma}{(n,a_1)\neq (0,0)}} \log |1-q_c(P)^{n+\tila_1}e^{2\pi i \tila_2}|
+
\sum_{\genfrac{}{}{0pt}{}{\bfa \in \OO\sigma}{}} \log |1-q_c(P)^{n+1-\tila_1}e^{-2\pi i \tila_2}|
\right)
\\& \hphantom{=} +
O_1\left( m \sum_{\bfa \in \OO\sigma}\frac{|q_c(P)|^{\nmax+\tila_1}+|q_c(P)|^{\nmax+1-\tila_1}}{(1-|q_c(P)|)^2}\right),
\end{aligned}
\end{equation*}
In addition to this, 
{\small
\begin{equation*}
%\label{eqDerivedLogUOofP}
\begin{aligned}
\left(\frac{d}{dq_c}\log |u|\right)(P)&= \frac{\ord_cu}{p} \re\frac{1}{q_c(P)}\ +
\\& \hphantom{=} +  
m\sum_{n=0}^{\nmax-1}
\sum_{\bfa \in \OO\sigma}
\re \left(
%\log |1-q_c^{n+\tila_1}e^{2\pi i \tila_2}|
\frac{-(n+\tila_1)q_c(P)^{n-1+\tila_1}e^{2\pi i\tila_2}}{1-q_c(P)^{n+\tila_1}e^{2\pi i \tila_2}}
+
%\log |1-q_c^{n+1-\tila_1}e^{-2\pi i \tila_2}|
\frac{-(n+1-\tila_1)q_c(P)^{n-\tila_1}e^{-2\pi i\tila_2}}{1-q_c(P)^{n+1-\tila_1}e^{-2\pi i \tila_2}}
\right)
\\& \hphantom{=} +
O_1\left( m
%\bigg(
%\frac{1}{(1-|q_c|)^2}-\sum_{n=1}^{\nmax-1}n|q_c|^{n-1}
%\bigg)
\frac{\nmax|q_c(P)|^{\nmax-1}-(\nmax-1)|q_c(P)|^{\nmax}}{(1-|q_c(P)|)^2}
\sum_{\bfa \in \OO\sigma}
\left(\frac{1}{(1-|q_c(P)|^{\nmax+\tila_1})}+\frac{1}{1-|q_c(P)|^{\nmax+1-\tila_1}}\right)
\right).
\end{aligned}
\end{equation*} 
}%
\end{corollary}

\begin{proof}
As in the proof of Proposition~\ref{passu} we may assume that ${\sigma=I}$ and~$c$ is the cusp represented by $i\infty$. In this case the result is an immediate consequence of Propositions~\ref{papprsiegel} and~\ref{pderivsieg}. 
\end{proof}

In the program, we use the above estimate  for $\log |u(P)|$ first
with ${\nmax = 2}$, and then with larger and
larger $\nmax$ when the overall precision needs to be increased. As for the estimate 
 for $\left(\frac{d}{dq_c}\log|u|\right)(P)$, we use it only with ${\nmax=1}$, which turns out to be sufficient whenever we used it.
For larger~$p$ one may want to take a slightly better approximation, for instance, with ${\nmax=2}$.

%%%%%%%%%%%%%%%%%%%%%%%%%%%%%%%%%%%%%%%%%%%%%%%%%%%%%%%%%%%%%%%%%

{\footnotesize
\bibliographystyle{plain}
\bibliography{nonsplit}

\begin{thebibliography}{10}

\bibitem{BS13}
Aur\'{e}lien Bajolet and Min Sha.
\newblock Bounding the {$j$}-invariant of integral points on
  {$X_{\mathrm{ns}}^+(p)$}.
\newblock {\em Proc. Amer. Math. Soc.}, 142(10):3395--3410, 2014.

\bibitem{BD69}
Alan Baker and Harold Davenport.
\newblock The equations {$3x^{2}-2=y^{2}$} and {$8x^{2}-7=z^{2}$}.
\newblock {\em Quart. J. Math. Oxford Ser. (2)}, 20:129--137, 1969.

\bibitem{BDMTV19}
Jennifer Balakrishnan, Netan Dogra, J.~Steffen M\"{u}ller, Jan Tuitman, and Jan
  Vonk.
\newblock Explicit {C}habauty-{K}im for the split {C}artan modular curve of
  level 13.
\newblock {\em Ann. of Math. (2)}, 189(3):885--944, 2019.

\bibitem{Ba09}
Burcu Baran.
\newblock A modular curve of level 9 and the class number one problem.
\newblock {\em J. Number Theory}, 129(3):715--728, 2009.

\bibitem{Ba10}
Burcu Baran.
\newblock Normalizers of non-split {C}artan subgroups, modular curves, and the
  class number one problem.
\newblock {\em J. Number Theory}, 130(12):2753--2772, 2010.

\bibitem{Ba14}
Burcu Baran.
\newblock An exceptional isomorphism between modular curves of level 13.
\newblock {\em J. Number Theory}, 145:273--300, 2014.

\bibitem{Bi95}
Yuri Bilu.
\newblock Effective analysis of integral points on algebraic curves.
\newblock {\em Israel J. Math.}, 90(1-3):235--252, 1995.

\bibitem{Bi02}
Yuri Bilu.
\newblock Baker's method and modular curves.
\newblock In {\em A panorama of number theory or the view from {B}aker's garden
  ({Z}\"{u}rich, 1999)}, pages 73--88. Cambridge Univ. Press, Cambridge, 2002.

\bibitem{BH96}
Yuri Bilu and Guillaume Hanrot.
\newblock Solving {T}hue equations of high degree.
\newblock {\em J. Number Theory}, 60(2):373--392, 1996.

\bibitem{BH98}
Yuri Bilu and Guillaume Hanrot.
\newblock Solving superelliptic {D}iophantine equations by {B}aker's method.
\newblock {\em Compositio Math.}, 112(3):273--312, 1998.

\bibitem{BH99}
Yuri Bilu and Guillaume Hanrot.
\newblock Thue equations with composite fields.
\newblock {\em Acta Arith.}, 88(4):311--326, 1999.

\bibitem{BHV01}
Yuri Bilu, Guillaume Hanrot, and Paul~M. Voutier.
\newblock Existence of primitive divisors of {L}ucas and {L}ehmer numbers.
\newblock {\em J. Reine Angew. Math.}, 539:75--122, 2001.
\newblock With an appendix by M. Mignotte.

\bibitem{BI11}
Yuri Bilu and Marco Illengo.
\newblock Effective {S}iegel's theorem for modular curves.
\newblock {\em Bull. Lond. Math. Soc.}, 43(4):673--688, 2011.

\bibitem{BP11}
Yuri Bilu and Pierre Parent.
\newblock Runge's method and modular curves.
\newblock {\em Int. Math. Res. Not. IMRN}, (9):1997--2027, 2011.

\bibitem{BP11a}
Yuri Bilu and Pierre Parent.
\newblock Serre's uniformity problem in the split {C}artan case.
\newblock {\em Ann. of Math. (2)}, 173(1):569--584, 2011.

\bibitem{BPR13}
Yuri Bilu, Pierre Parent, and Marusia Rebolledo.
\newblock Rational points on {$X^+_0(p^r)$}.
\newblock {\em Ann. Inst. Fourier (Grenoble)}, 63(3):957--984, 2013.

\bibitem{Br71}
Richard~P. Brent.
\newblock An algorithm with guaranteed convergence for finding a zero of a
  function.
\newblock {\em Comput. J.}, 14:422--425, 1971.

\bibitem{Ca19}
Yulin Cai.
\newblock An explicit bound for integral points on modular curves.
\newblock {\href{http://arxiv.org/abs/1910.10405}{arxiv:1910.10405}}, 2019.

\bibitem{CC04}
Imin Chen and Chris Cummins.
\newblock Elliptic curves with nonsplit mod 11 representations.
\newblock {\em Math. Comp.}, 73(246):869--880, 2004.

\bibitem{sage}
The~Sage Developers.
\newblock {\em {S}age{M}ath ({V}ersion 7.3)}, 2017.
\newblock \url{http://www.sagemath.org}.

\bibitem{DS05}
Fred Diamond and Jerry Shurman.
\newblock {\em A first course in modular forms}, volume 228 of {\em Graduate
  Texts in Mathematics}.
\newblock Springer-Verlag, New York, 2005.

\bibitem{FP85}
Ulrich Fincke and Michael~E. Pohst.
\newblock Improved methods for calculating vectors of short length in a
  lattice, including a complexity analysis.
\newblock {\em Math. Comp.}, 44(170):463--471, 1985.

\bibitem{Ha00}
Guillaume Hanrot.
\newblock Solving {T}hue equations without the full unit group.
\newblock {\em Math. Comp.}, 69(229):395--405, 2000.

\bibitem{He52}
Kurt Heegner.
\newblock Diophantische {A}nalysis und {M}odulfunktionen.
\newblock {\em Math. Z.}, 56:227--253, 1952.

\bibitem{vKM16}
Rafael~von K\"anel and Benjamin Matschke.
\newblock {S}olving {$S$}-unit, {M}ordell, {T}hue, {T}hue--{M}ahler and
  generalized {R}amanujan--{N}agell equations via {S}himura--{T}aniyama
  conjecture.
\newblock {\href{http://arxiv.org/abs/1605.06079}{arxiv:1605.06079}}, 2016.

\bibitem{Ke85}
Monsur~A. Kenku.
\newblock A note on the integral points of a modular curve of level {$7$}.
\newblock {\em Mathematika}, 32(1):45--48, 1985.

\bibitem{KS10}
Ja~Kyung Koo and Dong~Hwa Shin.
\newblock On some arithmetic properties of {S}iegel functions.
\newblock {\em Math. Z.}, 264(1):137--177, 2010.

\bibitem{KL75}
Dan Kubert and Serge Lang.
\newblock Units in the modular function field. {I}.
\newblock {\em Math. Ann.}, 218(1):67--96, 1975.

\bibitem{KL81}
Daniel~S. Kubert and Serge Lang.
\newblock {\em Modular units}, volume 244 of {\em Grundlehren der
  Mathematischen Wissenschaften}.
\newblock Springer-Verlag, New York-Berlin, 1981.

\bibitem{La87}
Serge Lang.
\newblock {\em Elliptic functions}, volume 112 of {\em Graduate Texts in
  Mathematics}.
\newblock Springer-Verlag, New York, second edition, 1987.
\newblock With an appendix by J. Tate.

\bibitem{La95}
Serge Lang.
\newblock {\em Introduction to modular forms}, volume 222 of {\em Grundlehren
  der Mathematischen Wissenschaften}.
\newblock Springer-Verlag, Berlin, 1995.
\newblock With appendixes by D. Zagier and Walter Feit, Corrected reprint of
  the 1976 original.

\bibitem{LN97}
Rudolf Lidl and Harald Niederreiter.
\newblock {\em Finite fields}, volume~20 of {\em Encyclopedia of Mathematics
  and its Applications}.
\newblock Cambridge University Press, Cambridge, second edition, 1997.
\newblock With a foreword by P. M. Cohn.

\bibitem{Mas78}
John~Myron Masley.
\newblock Class numbers of real cyclic number fields with small conductor.
\newblock {\em Compositio Math.}, 37(3):297--319, 1978.

\bibitem{Ma78}
Barry Mazur.
\newblock Rational isogenies of prime degree (with an appendix by {D}.
  {G}oldfeld).
\newblock {\em Invent. Math.}, 44(2):129--162, 1978.

\bibitem{Mi15}
John~C. Miller.
\newblock Real cyclotomic fields of prime conductor and their class numbers.
\newblock {\em Math. Comp.}, 84(295):2459--2469, 2015.

\bibitem{PS87}
Attila Peth\H{o} and Rolf Schulenberg.
\newblock Effektives {L}\"{o}sen von {T}hue {G}leichungen.
\newblock {\em Publ. Math. Debrecen}, 34(3-4):189--196, 1987.

\bibitem{ST12}
Ren\'{e} Schoof and Nikos Tzanakis.
\newblock Integral points of a modular curve of level 11.
\newblock {\em Acta Arith.}, 152(1):39--49, 2012.

\bibitem{Se97}
Jean-Pierre Serre.
\newblock {\em Lectures on the {M}ordell-{W}eil theorem}.
\newblock Aspects of Mathematics. Friedr. Vieweg \& Sohn, Braunschweig, third
  edition, 1997.
\newblock Translated from the French and edited by Martin Brown from notes by
  Michel Waldschmidt, With a foreword by Brown and Serre.

\bibitem{Sh13}
Min Sha.
\newblock {\em Topics in Elliptic and Modular Curves}.
\newblock PhD thesis, Universit{é} Bordeaux 1, 2013.

\bibitem{Sh14a}
Min Sha.
\newblock Bounding the {$j$}-invariant of integral points on certain modular
  curves.
\newblock {\em Int. J. Number Theory}, 10(6):1545--1551, 2014.

\bibitem{Sh14}
Min Sha.
\newblock Bounding the {$j$}-invariant of integral points on modular curves.
\newblock {\em Int. Math. Res. Not. IMRN}, (16):4492--4520, 2014.

\bibitem{Shi71}
Goro Shimura.
\newblock {\em Introduction to the arithmetic theory of automorphic functions}.
\newblock Publications of the Mathematical Society of Japan, No. 11. Iwanami
  Shoten, Publishers, Tokyo; Princeton University Press, Princeton, N.J., 1971.
\newblock Kan\^{o} Memorial Lectures, No. 1.

\bibitem{Si68}
Carl~Ludwig Siegel.
\newblock Zum {B}eweise des {S}tarkschen {S}atzes.
\newblock {\em Invent. Math.}, 5:180--191, 1968.

\bibitem{pari}
{The PARI Group}.
\newblock {\em {PARI/GP} ({V}ersion 2.9.2)}, 2017.
\newblock Bordeaux, \url{http://pari.math.u-bordeaux.fr/}.

\bibitem{TW89}
Nikos Tzanakis and Benjamin M.~M. de~Weger.
\newblock On the practical solution of the {T}hue equation.
\newblock {\em J. Number Theory}, 31(2):99--132, 1989.

\bibitem{Wa97}
Lawrence~C. Washington.
\newblock {\em Introduction to cyclotomic fields}, volume~83 of {\em Graduate
  Texts in Mathematics}.
\newblock Springer-Verlag, New York, second edition, 1997.

\bibitem{Yu07}
Kunrui Yu.
\newblock {$p$}-adic logarithmic forms and group varieties. {III}.
\newblock {\em Forum Math.}, 19(2):187--280, 2007.

\end{thebibliography}
}

\vspace{2em}

\noindent
\begin{minipage}{12em}
Aur\'e{}lien Bajolet\\
{\addresssize\textit{Lyc\'e{}e Jean Dautet \\ 18 rue Delayant \\ 17000 La Rochelle, France\\
\hphantom{17000 La Rochelle, France }} }
\end{minipage}
\begin{minipage}{17em}
Yuri Bilu\\
{\addresssize\textit{Institut de Math\'e{}matiques de Bordeaux \\ Universit\'e{} de Bordeaux et CNRS \\ 351 cours de la Lib\'e{}ration \\ 33405 Talence, France}}
\end{minipage}
\begin{minipage}{17em}
Benjamin Matschke\\
{\addresssize\textit{Department of Mathematics and Statistics \\ Boston University \\ 111 Cummington Mall \\ Boston, MA 02215, USA}}
\end{minipage}

\end{document}